\documentclass{article}
\usepackage{graphicx} 
\usepackage{appendix}

\usepackage{pstricks,amsmath,amsfonts,latexsym,enumerate}
\usepackage[hyphens]{url}
\usepackage{epsfig,subfig,textcomp}
\usepackage{amsthm,subfig,longtable,framed,amsbsy,amssymb,amscd}
\usepackage{bbm}
\usepackage{multirow}
\usepackage{epstopdf}
\usepackage[T1]{fontenc}
\usepackage{palatino,cite}

\usepackage{amssymb,amscd}
\usepackage{nomencl}
\makenomenclature
\makeatletter
\def\myVCENTER#1{\vcenter{\hbox{$\m@th#1$}}}
\makeatother

\long\def\symbolfootnote[#1]#2{\begingroup\def\thefootnote{\fnsymbol{footnote}}\footnote[#1]{#2}\endgroup}

\definecolor{shadecolor}{gray}{0.99}
{\endMakeFramed}

\long\def\symbolfootnote[#1]#2{\begingroup\def\thefootnote{\fnsymbol{footnote}}\footnote[#1]{#2}\endgroup}
\def\qed{\hfill{$\vcenter{\hrule height1pt \hbox{\vrule width1pt height5pt
    \kern5pt \vrule width1pt} \hrule height1pt}$} \medskip}
\newcommand{\m}[1]{{\bf{#1}}}

\newcommand{\C}[1]{{\cal {#1}}}

\title{\bf Method for Constrained Computational Guidance and Control with Application to Hypersonic Entry}

\author{Emily M. Palmer\footnote{PhD Candidate, Department of Mechanical and Aerospace Engineering.  E-mail: emilypalmer@ufl.edu.}~ and Anil V. Rao\footnote{Professor, Department of Aerospace Engineering.  Fellow, American Astronautical Society.   Associate Fellow AIAA. \newline\hspace*{12pt} E-mail: anilvrao@ufl.edu.} \vspace{12pt} \\ {\em Department of Mechanical and Aerospace Engineering} \\ {\em University of Florida} \\ {\em Gainesville, FL 32611-6250}}

\date{}

\begin{document}

\maketitle

\begin{abstract}
A method for constrained computation guidance and control is described. The method employs guidance updates by periodically re-solving a constrained optimal control problem. Feasibility of the optimal control problem is improved by preventing constraint violations in the perturbed system. Constraint violation is prevented by augmenting the optimal control problem with a penalty term intended to increase the margin between a constraint and its limit when the perturbed system is within certain proximity of the constraint limit. The additional penalty term employs a smooth maximum approximation to increase the margin between the constraint and its limit. The method developed in this paper is demonstrated on the maximum cross-range reusable launch vehicle (RLV) entry problem, for which the reference solution contains a long-duration constrained arc. The results demonstrate that the method developed in this paper is viable for path-constrained optimal guidance and control. 
\end{abstract}

\renewcommand{\baselinestretch}{1}
\normalsize\normalfont

\section{Introduction}\label{section:introduction}

 Guidance is an on-board (real-time) process that computes commanded inputs to a dynamical system in order to direct the system towards a particular goal. The control system utilizes these guidance commands when executing the physical maneuvers using on-board hardware.  Due to limited computational resources, early guidance algorithms employed closed-form guidance laws that executed simple algebraic operations. These algorithms computed commanded inputs at low rates and could be executed using low-speed on-board processors.  Key examples of early guidance systems include U.S.~Space Shuttle ascent \cite{Jaggers74} and entry \cite{Harpold1981}.  In recent years, the increased availability of computational resources has made it possible to develop advanced guidance algorithms, known as {\em computational guidance and control} (CG\&C)\cite{Lu2017a,tsiotras2017toward}.  Specifically, CG\&C methods employ modern numerical methods to improve system performance. These methods are developed around the model/data itself, and are less specialized for use with a particular vehicle design. Additionally, these numerical methods allow a wider range of more complex problems to be solved as compared to analytical-based methods. 
 
 Advanced CG\&C methods are especially suitable for the complex problems that arise in aerospace applications. In particular, many aerospace systems have highly nonlinear dynamics, are subject to critical path constraints, and must operate in a wide range of environmental conditions.  One important application of CG\&C is atmospheric entry. Early guidance algorithms for entry were based on simple models and employed analytical expressions. For example, the US space shuttle guidance algorithm \cite{Harpold1981,lu2014entry} linked analytical deceleration functions of the Earth-relative velocity and the energy with respect to the Earth, based upon the Mach number, to create a reference drag profile.  Minor adjustments could then be made to attain the required downrange. With the development of advanced CG\&C methods, more sophisticated guidance algorithms for atmospheric entry can be developed that employ numerical methods to generate trajectories that could potentially optimize performance while satisfying critical constraints during entry.  

 To meet the needs for high-performance, in recent years optimization-based CG\&C methods have been developed. Many of the recent approaches reduce the computational complexity of the optimization problem in some manner. Examples of such approaches include convex optimization-based guidance  \cite{wang2024survey,liu2017survey,foust2020optimal} and optimal control-based methods, including model predictive control (MPC) \cite{qin2003survey,mayne2014,eren2017model} and model predictive static programming (MPSP) \cite{sachan2019waypoint,dwivedi2011suboptimal,kumar2019model}. Convex optimization-based guidance approaches transform nonconvex problems into convex optimization problems, which are guaranteed to converge in a finite number of steps for feasible problems. The convex relaxation of the problem may or may not be lossless \cite{accikmecse2011lossless,harris2014lossless} (that is, the solution to the convex optimization problem is the same as the solution of the original nonconvex problem). Additionally, if nonconvex path constraints are present and the convex relaxation of the path constraints is not exact, the solution may not be feasible with respect to the original constraints. MPC approaches reduce the computational burden of the optimization problem by restricting the prediction horizon and linearization of the problem. Recent efforts have been made towards developing nonlinear MPC methods \cite{allgower2012nonlinear}, though the application to complex systems is limited and the incorporation of nonlinear inequality constraints can be challenging. MPSP approaches solve a static optimization problem only in control variables, allowing for a closed-form solution of the suboptimal control history. MPSP requires the problem be linearized and is limited to applications where the objective is in terms of the control only.

As described above, one particular class of optimization-based CG\&C methods is based on solving an optimal control problem in real time.  Because most optimal control problems do not have analytic solutions, any high-performance CG\&C guidance algorithm based on optimal control must be numerical.  Numerical methods for optimal control fall into two broad categories: indirect or direct methods. In an indirect method, the necessary optimality conditions are derived using the calculus of variations and the optimal control problem is converted into a Hamiltonian boundary-value problem (HBVP), which is then solved numerically. It is noted, however, that for many problems it is challenging to derive these necessary conditions.  In addition, even if these necessary conditions can be derived, most indirect methods require a good initial guess (including a guess of the non-intuitive costate).  In a direct method, the state and control are approximated at a discrete set of points and the optimal control problem is transcribed into a finite-dimensional nonlinear programming problem (NLP). While direct methods often have a larger radius of convergence, some direct methods produce inaccurate solutions while others are not well-suited to problems with constraints. In addition, depending upon the optimization techniques employed, some direct methods may not be computationally efficient enough for use in on-board computation.  

Of the optimal control-based CG\&C methods, guidance methods that employ direct collocation have been proposed for nonlinear optimal control problems  \cite{bollino2006optimal,bollino2007pseudospectral,Dennis2019,dennis2021optimal}. These methods take advantage of the advancements that have been made in nonlinear optimization. Recently, the CG\&C method of Ref.~\cite{Dennis2019} was developed in which the optimal control problem is transcribed into a large sparse NLP using collocation at Legendre-Gauss-Radau (LGR) points.  The LGR transcription of the optimal control problem is then re-solved at each guidance update on the entire remaining horizon using the full nonlinear dynamics.  The method presented in Ref.~\cite{Dennis2019} has the key feature that it is able to quickly re-solve general nonconvex NLPs via a mesh truncation and remapping strategy, where the portion of the mesh associated with the expired horizon is deleted at the start of each guidance cycle.  While the method of Ref.~\cite{Dennis2019} employs the full nonlinear dynamic model, it is not designed to remain feasible with respect to inequality path constraints.  More specifically, when employing the method of Ref.~\cite{Dennis2019}, it is possible that path constraints may be violated when the actual system deviates from the reference optimal solution, and the optimal control problem being solved may become infeasible. A demonstration of this phenomenon is exhibited in the constrained Mars entry application given in Ref.~\cite{Palmerjan2022}.

In this paper, a new method is developed that is intended to prevent constraint violations by augmenting the objective functional of the optimal control problem with a penalty term intended to increase the margin between the constraint function and its limit. In the case of an upper limit, the penalty term is formulated such that the maximum value of the constraint function is minimized. Conversely, in the case of a lower limit, the penalty term is formulated to minimize the maximum value of the negative of the constraint function. To make this penalty term tractable with gradient-based optimization methods, the maximum function is replaced by a smooth approximation.  In particular, a form of the log-sum-exponential function, which has primarily been limited to use in machine learning applications \cite{blanchard2021accurately}, is employed. In order to avoid unnecessarily penalizing the problem when constraint violation is unlikely, this paper presents a method in which the weight on the penalty term in the objective functional is chosen at the start of each guidance cycle based on the proximity of the actual system to the constraint  limit. When the actual system is far from the constraint limit, the weight is chosen to vanish such that the reference performance index is prioritized. Conversely, when the actual system approaches the constraint limit, the weight on the penalty is increased to ensure a sufficient margin between the constraint function and its limit.  It is also noted that the method of this paper employs a similar mesh truncation and re-mapping strategy as that described in Ref.~\cite{Dennis2019}.  It is noted that highly preliminary work on the constrained optimal guidance method presented in this paper has been presented in Refs.~\cite{PalmerRao2023,palmer2022robust,palmer2025}. 

The contributions of this paper are as follows. First, a constrained optimization based CG\&C method is described that allows a nonconvex optimal control problem to be solved via transcription to an NLP using LGR collocation. The presented method allows for the direct incorporation of path constraints and boundary conditions as hard constraints. The NLP is solved over the entire remaining horizon at the start of each guidance cycle, and the computational burden of this NLP is reduced using the mesh truncation and re-mapping strategy of Ref.~\cite{Dennis2019}. The control obtained from the solution of the NLP is used to simulate the actual dynamical system (which is assumed to differ from the reference dynamical model), and the state of the actual system at the end of the guidance cycle is used to update the initial conditions of the NLP for the ensuing guidance cycle. Second, a new constraint penalty is described to prevent constraint violations within the actual system by increasing the margin between the constraint function and its limit to ensure safe operation and feasibility of the NLP at each guidance cycle. Third, this paper presents a strategy in which the weights in the modified objective functional are determined at the start of each guidance cycle based on the proximity of the actual system to the constraint limits. Choosing the weights at the beginning of each guidance cycle based on the needs of the actual system allows for a balance between maximizing performance and ensuring feasibility. The method described in this paper is then applied to a complex maximum cross-range reusable launch vehicle entry problem \cite{Betts2020}, for which the reference optimal solution contains a long-duration constrained arc, and the performance is compared against the method of Ref.~\cite{Dennis2019}. The method of this paper is first employed using an objective functional comprised of a weighted combination of the reference objective and the constraint function penalty, and the effect of varying these weights is studied. Then, an additional penalty term is introduced into the objective functional, intended to reduce oscillatory behavior that occurs as a result of incorporation of the constraint function penalty term. Finally, an extensive Monte-Carlo simulation shows that the method of this paper eliminates infeasibility due to path constraint violation, while the method described in Ref.~\cite{Dennis2019} often suffers from constraint violation. 

This paper is organized as follows. Section~\ref{sect:background} provides the background material for the development of the method described in this paper.  Section~\ref{sect:previous-cgs-method} provides a summary of the computational optimal guidance algorithm developed in Ref.~\cite{Dennis2019} that provides a starting point for the method of this paper.  Section~\ref{sect:new-cgs-method} describes the new constrained CG\&C method of this paper. Section~\ref{Example Application} describes the atmospheric entry problem which is used as the basis for demonstrating the effectiveness of the method developed in Section~\ref{sect:new-cgs-method}.  Section~\ref{results} provides the details of the analysis of the method of Section~\ref{sect:new-cgs-method} on the atmospheric entry problem described in Section~\ref{Example Application}.  Section~\ref{discussion} provides a discussion on the results given in Section~\ref{results}. Section \ref{conclusions} provides conclusions on this research.

\section{Background Material}\label{sect:background}

This section provides the background material for the method described in Section \ref{sect:new-cgs-method}. This section is divided into the following three subsections.  First, notation and conventions are described.  Second, the dynamic models used in the development of the method are described.  Third, the general constrained optimal control problem is described.  

\subsection{Notation and Conventions}

The following notation and conventions are employed throughout this paper. To begin, the independent variable $t\in[t_0,t_f]$ denotes the elapsed time, where $t_0$ is the start (initial) time and $t_f$ is the terminal time.  Next, the dynamic system is assumed to be described by a system of first-order differential equations
\begin{equation}\label{ode}
    \frac{d\m{x}(t)}{dt} = \dot{\m{x}}(t) = \m{f}(\m{x}(t),\m{u}(t),t), 
\end{equation}
where $\m{x}(t)\in\mathbb{R}^{n_x}$ is the state and $\m{u}(t)\in\mathbb{R}^{n_u}$ is the control.  At any instant of time, the state and control are assumed to be column vectors, respectively, of the form
\begin{equation}
\begin{array}{lcl}
    \m{x}(t)=\left[\begin{array}{c}x_1(t) \\ \vdots \\ x_{n_x}(t)\end{array}\right]  
& , &    \m{u}(t)=\left[\begin{array}{c}u_1(t) \\ \vdots \\ u_{n_u}(t)\end{array}\right].
\end{array}
\end{equation}
Next, the vector field $\m{f}:\mathbb{R}^{n_x}\times\mathbb{R}^{n_u}\times\mathbb{R}\rightarrow\mathbb{R}^{n_x}$ that defines the differential equations given in Eq.~\eqref{ode} has the form
\begin{equation}\label{vector-field-f}
    \m{f}(\m{x}(t),\m{u}(t),t) = 
    \left[
    \begin{array}{c}
        f_1(\m{x}(t),\m{u}(t),t) \\ \vdots \\ f_{n_x}(\m{x}(t),\m{u}(t),t)
    \end{array}
    \right].
\end{equation}
In addition, if a system is subject to inequality path constraints of the form
\begin{equation}
    \m{c}_{\min}\leq\m{c}(\m{x}(t),\m{u}(t),t)\leq\m{c}_{\max},
\end{equation}
then the function $\m{c}:\mathbb{R}^{n_x}\times\mathbb{R}^{n_u}\times\mathbb{R}\rightarrow\mathbb{R}^{n_c}$ that defines these path constraints has the form
\begin{equation}\label{vector-field-c}
    \m{c}(\m{x}(t),\m{u}(t),t) = 
    \left[
    \begin{array}{c}
        c_1(\m{x}(t),\m{u}(t),t) \\ \vdots \\ c_{n_c}(\m{x}(t),\m{u}(t),t)
    \end{array}
    \right]. 
\end{equation}

\subsection{Dynamic Models}

Two dynamic models will be employed in this paper.  The first dynamic model is referred to as the {\em reference dynamic model} and has the form given in Eq.~\eqref{ode}. For a given control, $\m{u}(t)$, the solution of Eq.~\eqref{ode} is referred to as the {\em reference state}, $\m{x}(t)$.  The second dynamic model is referred to as the {\em perturbed dynamic model} and has the form
\begin{equation}\label{ode-perturbed}
  \frac{d\tilde{\m{x}}}{dt}= \dot{\tilde{\m{x}}}(t) =\tilde{\m{f}}(\tilde{\m{x}}(t),\m{u}(t),t).
\end{equation}
For a given control, $\m{u}(t)$, the solution of Eq.~\eqref{ode-perturbed} is referred to as the {\em perturbed state}, $\tilde{\m{x}}(t)$.  As will be described in Section \ref{sect:ocp}, the reference dynamic model given in Eq.~\eqref{ode} is used to solve an optimal control problem that forms the basis of the guidance method of this paper, while the perturbed dynamic model given in Eq.~\eqref{ode-perturbed} is the model that is used to simulate the motion of the "actual" system. Note that the reference and perturbed dynamic models differ due to environmental disturbances, measurement errors, or other modeling imperfections.   

\subsection{Constrained Optimal Control Problem}\label{sect:ocp}

The basis of the constrained optimal guidance method of this paper is the solution of a constrained nonlinear optimal control problem.  To this end, consider the following general optimal control problem in Bolza form. Minimize the objective functional 
\begin{equation}\label{eq:obj}
\C{J} = \C{M}(\m{x}(t_f),t_f) + \int_{t_0}^{t_f} \C{L}(\m{x}(t),\m{u}(t),t)\, dt,
\end{equation}
subject to the dynamic constraints (which is a rearrangement of Eq.~\eqref{ode}),
\begin{equation}\label{eq:dyncon}
  \frac{d\m{x}}{dt}-\m{f}(\m{x}(t),\m{u}(t),t)=\m{0}, 
\end{equation}
the inequality path constraints
\begin{equation}\label{eq:pathcon}
   \m{c}_{\textrm{min}} \leq \m{c}(\m{x}(t),\m{u}(t),t) \leq \m{c}_{\textrm{max}}, \; 
\end{equation}
and the boundary conditions 
\begin{equation}\label{eq:boundcond}
\begin{array}{rcl}
\m{x}(t_0) &=& \m{x}_0, \\
 \m{b}(\m{x}(t_f),t_f) &=& \m{0}, \\
t_0 &= &\mathrm{Fixed}.
\end{array}
\end{equation}
The function $\C{M}: \mathbb{R}^{n_x} \times \mathbb{R} \rightarrow \mathbb{R}$ is the Mayer (terminal) cost, $\C{L}: \mathbb{R}^{n_x} \times \mathbb{R}^{n_u} \times \mathbb{R} \rightarrow \mathbb{R}$ is the Lagrangian (running cost), $\m{b} \in \mathbb{R}^{n_x} \times \mathbb{R} \rightarrow \mathbb{R}^{n_b}$ is the function that describes the $n_b$ terminal constraints, $\m{x}_0\in\mathbb{R}^{n_x}$ is a specified initial state, and $t_0\in\mathbb{R}$ is a specified initial time.  The solution to the optimal control problem given in Eqs. \eqref{eq:obj}-\eqref{eq:boundcond} is defined to be $\C{B} = (\m{x}^*(t),\m{u}^*(t),t_f^*)$, where $\m{x}^*(t)$ is the {\em optimal state}, and $\m{u}^*(t)$ is the {\em optimal control}. The optimal state $\m{x}^*(t)$ evolves according to reference dynamics $\m{f}(\m{x}(t),\m{u}(t),t)$ given in Eq.~\eqref{ode} and is referred to as the {\em reference state}. 


\section{Review of CG\&C Method of Ref.~\cite{Dennis2019}}\label{sect:previous-cgs-method}

This section provides a summary of the computational guidance and control method  described in Ref.~\cite{Dennis2019} for problems without path constraints, that is, the inequality path constraints defined by Eq.~\eqref{eq:pathcon} are {\em not} included.  In the case where no path constraints are included, the method of Ref.~\cite{Dennis2019} solves an optimal control problem where it is desired to minimize the objective functional in Eq.~\eqref{eq:obj} subject to the dynamic constraints of Eq.~\eqref{eq:dyncon} and the boundary conditions of Eq.~\eqref{eq:boundcond}.  This path-unconstrained optimal control problem is solved every $\Delta T$ time units, where the duration $\Delta T$ is referred to as a {\em guidance cycle}. Re-solving the optimal control problem at these specified update times is referred to a {\em guidance update} with start and terminus times given, respectively, as
\begin{eqnarray}\label{guidance-update-start-end}
    t_0^{(s)} & = & t_0 + s \Delta T, \\
    t_e^{(s)} & = & t_0 + (s+1)\Delta T. 
\end{eqnarray}
where $s\in[1,2,\ldots,S]$ is the guidance update number, and $S$ is the total number of guidance updates.

Suppose that $\C{B}_u^{(s)}=(\m{x}^*(t),\m{u}^*(t),{t_f^*}^{(s)})$ is the solution of the path-unconstrained optimal control problem defined by Eqs.~\eqref{eq:obj}, \eqref{eq:dyncon}, and \eqref{eq:boundcond} on the horizon $t\in [t_0^{(s)},{t_f^*}^{(s)}]$. In the method of Ref.~\cite{Dennis2019}, the optimal control problem is solved using LGR collocation \cite{Garg2011a,Garg2011b,Garg2010,Darby2010,Darby2011b,Patterson2015,Liu2015,Liu2018}.  The initial mesh is solved to meet all required accuracy tolerances so as to eliminate any unnecessary mesh refinement when re-solving the optimal control problem.  The optimal control $\m{u}^*(t)$ is applied to the actual system over guidance cycle $s$, which evolves according to the perturbed dynamics given by Eq.~\eqref{ode-perturbed}. At the terminus time of guidance cycle $s$, the mesh used to obtain the solution $B_u^{(s)}$ is truncated and re-mapped to the unexpired horizon $[t_e^{(s)},{t_f^*}^{(s)}]$. The truncated mesh is used, together with the reference state and optimal control on the unexpired horizon, as an initial guess to re-solve the optimal control problem over the time horizon $[t_0^{(s+1)},{t_f}^{(s+1)}]$, with the updated initial condition $\m{x}_0^{(s+1)} = \tilde{\m{x}}(t_e^{(s)})$. The details of the mesh truncation and re-mapping strategy are outside of the scope of this paper but a full description can be found in Ref.~\cite{Dennis2019}. This process of re-solving the optimal control problem at the start of each guidance cycle, implementing the resulting optimal control using the perturbed dynamics over the ensuing guidance cycle, and truncating and re-mapping the mesh is repeated for each guidance update $s\in [1,2,\ldots,S]$.  


\section{CCG\&C Method for Constrained Optimal Guidance and Control}\label{sect:new-cgs-method}

The method of Ref.~\cite{Dennis2019} needs to be modified significantly when considering applications with path constraints.  The reason for this significant modification is that the solution of a path-constrained optimal control problem $\C{B}^{(s)}$ may contain segments that lie along a constrained arc where a path constraint is active for a nonzero duration of time.  Note that near or on the constrained arc, maintaining feasibility in the perturbed system with respect to path constraints can be challenging, if not impossible, because disturbances or modeling errors may lead to violations of the path constraints. At times where the path constraint is violated, the updated initial condition $\m{x}_0^{(s+1)} = \tilde{\m{x}}(t_e^{(s)})$ will result in an infeasible optimal control problem.  This section describes a new method for constrained computational guidance and control (hereby referred to as the CCG\&C method) that is designed to prevent constraint violation and maintain feasibility across the entire horizon, and this new method is inspired by the work of Ref.~\cite{Dennis2019}.

Suppose that the optimal control problem defined in Section \ref{sect:ocp} includes inequality path constraints and is now defined by Eqs.~\eqref{eq:obj}--\eqref{eq:boundcond}.  Furthermore, suppose that is it desired to maintain feasibility with respect to these path constraints across the entire horizon.  Consider the case where the solution $\C{B}^{(s)}$ contains a segment $t\in[t_1,t_2]$ such that one of the path constraints $c_i,\;i\in[1,\ldots,n_c]$ is active, that is, the solution $\C{B}^{(s)}$ contains a constrained arc such that either
\begin{displaymath}
c_i(\m{x}^*(t),\m{u}^*(t),t) = c_{i,\max},
\end{displaymath}
or 
\begin{displaymath}
c_i(\m{x}^*(t),\m{u}^*(t),t) = c_{i,\min}.
\end{displaymath}
If the dynamic system is modeled perfectly and is subject to no other errors (for example, measurement errors), theoretically the reference state, $\m{x}^*(t)$, and the perturbed state, $\tilde{\m{x}}(t)$, will be identical.  In reality, however, $\tilde{\m{x}}(t)$ will differ from $\m{x}^*(t)$, and the actual system may violate the path constraint defined by the function $c_i,\;i\in[1,\ldots,n_c]$. As a result, for the perturbed solution  $\tilde{\m{x}}(t)$ on the interval $t\in[t_1,t_2]$, it can be the case that either
\begin{displaymath}
  c_i(\tilde{\m{x}}(t),\m{u}^*(t),t) > c_{i,\max},
\end{displaymath}
or
\begin{displaymath}
  c_i(\tilde{\m{x}}(t),\m{u}^*(t),t) < c_{i,\min}. 
\end{displaymath}
Consequently, on a segment where the perturbed state violates the constraints, the initial condition at the start of some guidance cycle $s$, that is, $\m{x}_0^{(s)}=\tilde{\m{x}}(t_e^{(s-1)})$, will be such that the optimal control problem defined in Section \ref{sect:ocp} will be infeasible.  

The approach developed in this paper for avoiding infeasibilities with respect to the path constraints is to modify the optimal control problem in real-time to prevent the perturbed state from violating the path constraints.  Specifically, the objective functional of the optimal control problem defined in Section \ref{sect:ocp} is modified to increase the margin (in the feasible direction) between the path constraint function and the path constraint limit.  This modified objective functional, $\C{J}_a$, has the form
\begin{equation}\label{augmented-objective-ocp}
    \C{J}_a = w_1^{(s)} \C{J} + w_2^{(s)} \C{J}_p,
\end{equation}
where $\C{J}$ is the reference objective functional, $\C{J}_p$ is a penalty term, and $(w_1^{(s)},w_2^{(s)})$ are weights such that $|w_i^{(s)}|\leq 1$ and 
\begin{equation}\label{eq:sum-of-weights}
    w_1^{(s)} + w_2^{(s)} = 1. 
\end{equation}
It is noted that the weights $(w_1^{(s)},w_2^{(s)})$ need not be constants but can be functions of the value of the constraint function, evaluated at the terminus of the previous guidance cycle using the perturbed state.  In this paper, the weight $w_2^{(s)}$ is chosen such that its value increases as the distance between the path constraint function and its limit decreases.  Then, because of the relationship between $w_1^{(s)}$ and $w_2^{(s)}$ given in Eq.~\eqref{eq:sum-of-weights}, as $w_2^{(s)}$ increases, $w_1^{(s)}$ decreases (and vice versa). Determining the weights in real-time based on the current proximity of a path constraint function to its limit prioritizes the reference objective functional when the system is far from a constraint limit, while prioritizing the additional term when the path constraint function is in close proximity to its limit. In Section~\ref{subsec:comparing weights}, a threshold function will be utilized and defined in further detail to determine the objective functional weights.
The approach for choosing the additional penalty term $\C{J}_p$ in Eq.~\eqref{augmented-objective-ocp} is now described.  

\subsection{Penalty Term in Objective Functional}

To start, consider the following penalty term $\C{J}_p$, assuming that at the start of the current guidance cycle, $s$, the perturbed state is feasible.  If the path constraint function is approaching its upper limit, then the goal is to minimize the maximum value of the path constraint function
\begin{equation}\label{approaching-upper-limit}
\C{J}_p = \max_t \left\{  c_i \right\}.
\end{equation}
On the other hand, if the path constraint function is approaching its lower limit, then the goal is to maximize the minimum value of the path constraint function
\begin{equation}\label{approaching-lower-limit-pre}
\C{J}_p = \min_t \left\{ c_i \right\}. 
\end{equation}
Now, maximizing the minimum value of the path constraint function is equivalent to minimizing the maximum of the negative of the path constraint function in Eq.~\eqref{approaching-lower-limit-pre}. As such Eq.~\eqref{approaching-lower-limit-pre} can be re-written such that the goal is to minimize
\begin{equation}\label{approaching-lower-limit}
\C{J}_p = \max_t \left\{ -c_i \right\}. 
\end{equation}
Unifying each of the cases shown in Eqs.~\eqref{approaching-upper-limit} and \eqref{approaching-lower-limit} gives
\begin{equation}\label{unified-penalty}
    \C{J}_p = \left\{\begin{array}{lcl} \max_t \left\{ + c_i \right\} & , & \textrm{approaching upper limit}, \\ \max_t \left\{ -c_i \right\} & , & \textrm{approaching lower limit}. \end{array} \right.
\end{equation}
Another way of looking at Eq.~\eqref{unified-penalty} is that the goal of adding the penalty term $\C{J}_p$ is to maximize the minimum distance between the path constraint function and its limit.

It is important to see that the function in Eq.~\eqref{unified-penalty} is not differentiable because taking the maximum value over $t$ is itself a nonsmooth function.  As a result, it is not possible to use a gradient-based optimization method if the penalty term $\C{J}_p$ is added to the original objective functional $\C{J}$; therefore, an alternative approach must be developed.  In the next section the maximum function shown in Eq.~\eqref{unified-penalty} is replaced with a smooth approximation that will enable the inclusion of the penalty term in the objective functional of the optimal control problem. 

\subsection{Smooth Approximation of Nonsmooth Penalty Term in Objective Functional}

\subsubsection{Log-Integral-Exponential Function}

In order to employ a gradient-based optimization method, all functions must be differentiable.  Because the maximum function in Eq.~\eqref{unified-penalty} is not differentiable, it must be replaced by a smooth approximation.  Several smooth approximations to the maximum function have been developed. One particular smooth approximation of the maximum function that is well-behaved is the {\em log-sum-exponential approximation} \cite{blanchard2021accurately}. Note, however, that the log-sum-exponential approximation is formulated for discrete data as opposed to continuous functions. Because in this paper the optimal control problem is presented in its continuous form, it is necessary to modify the log-sum-exponential approximation for continuous functions.  In this paper, a particular variation of the log-sum-exponential function, referred to as the {\em log-integral-exponential} function and defined using a function  $P_\beta$, is utilized
\begin{equation}\label{log-integral-exponential-scaled}
P_{\beta}(y(t)) =\frac{1}{\beta} \log \left( \int_{t_0}^{t_f} \exp(\beta y(t)) dt  \right),
\end{equation}
where $\beta$ is a tuning parameter, and the approximation of the maximum value given by Eq.~\eqref{log-integral-exponential-scaled} increases in accuracy as $\beta$ increases.

\subsubsection{Smooth Penalty Term}

\begin{figure}[h]
\centering
\includegraphics[scale=0.33]{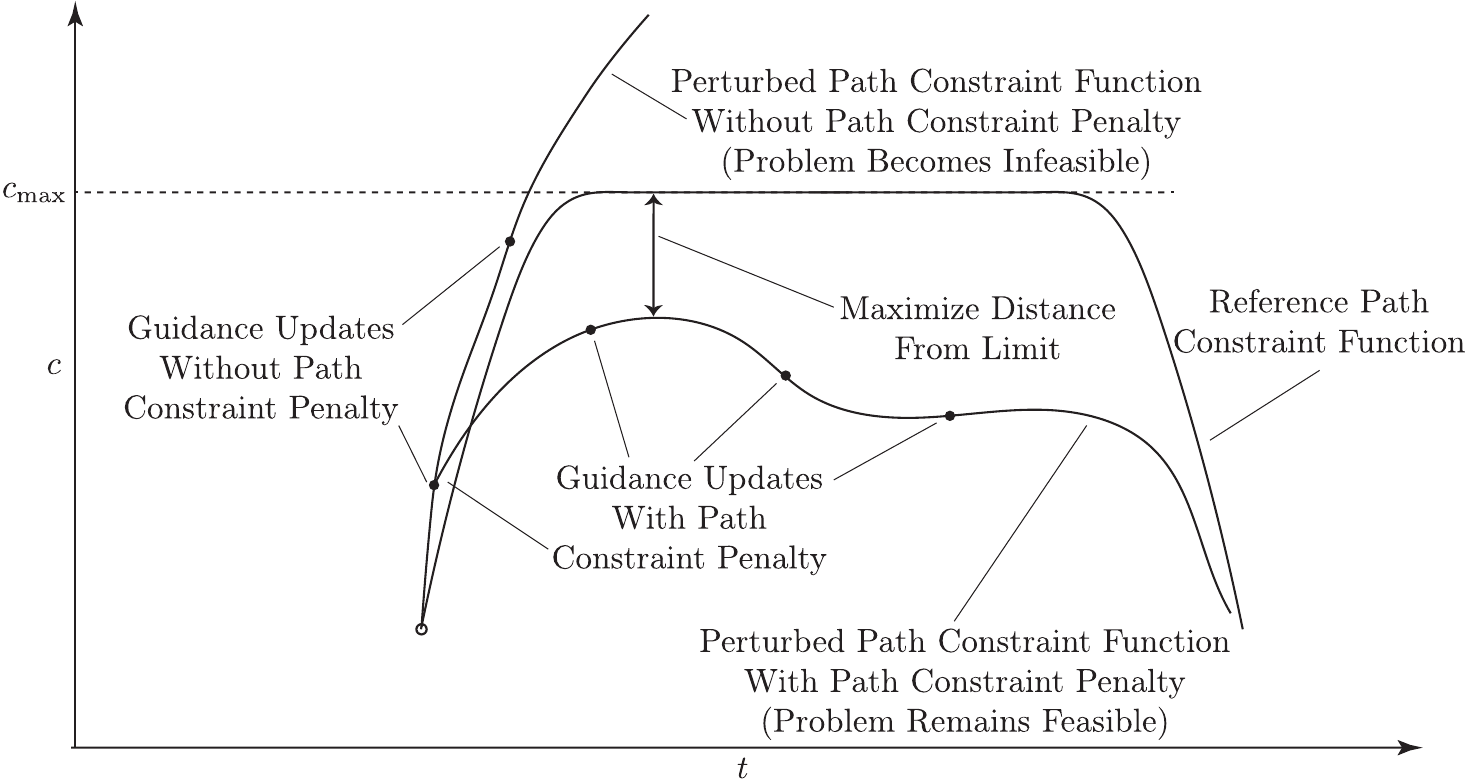}
\caption{Influence of penalty term in the objective functional of the constrained optimal control problem\label{fig:penalty-diagram}}
\end{figure}

In order to increase computational tractability and balance the terms in the modified objective functional, it is desirable to replace the penalty of Eq.~\eqref{unified-penalty} with a scaled path constraint function, such that
\begin{equation}\label{unified-penalty-scaled}
    \C{J}_p = \left\{\begin{array}{lcl} \max_t \left\{ + \frac{c_i}{\eta_i} \right\} & , & \textrm{approaching upper limit}, \\ \max_t \left\{ -\frac{c_i}{\eta_i} \right\} & , & \textrm{approaching lower limit}, \end{array} \right.
\end{equation}
where $\eta_i>0$ is a scaling parameter. The scaling parameter can be chosen to be the maximum absolute value of $c_i$ such that $|c_i(\m{x}(t),\m{u}(t),t)|/\eta_i \leq 1$.The penalty given by Eq.~\eqref{unified-penalty-scaled} is then replaced with the smooth maximum approximation given by Eq.~\eqref{log-integral-exponential-scaled}, where $y(t) = c_i/\eta_i$, such that
\begin{equation}\label{eq: scaled penalty term}
    \C{J}_{P_{\beta}} = \left\{\begin{array}{lcl} P_{\beta} \left(+\frac{c_i}{\eta_i}\right) & , & \textrm{approaching upper limit}, \\ P_{\beta} \left(-\frac{c_i}{\eta_i}\right) & , & \textrm{approaching lower limit}. \end{array} \right.
\end{equation}
 The function $\C{J}_{P_{\beta}}$ is used as the additional term in the objective functional given by Eq.~\eqref{augmented-objective-ocp}.  Figure~\ref{fig:penalty-diagram} provides a schematic of the evolution of the path constraint function in the absence and presence of the penalty term in Eq.~\eqref{eq: scaled penalty term} for the case of an upper limit. Note, Eq.~\eqref{eq: scaled penalty term} is defined for a single component $c_i(\m{x}(t),\m{u}(t),t),\; i\in[1,\ldots,n_c]$, of the path constraint function $\m{c}(\m{x}(t),\m{u}(t),t)$.   If two or more components of the path constraint function need to be penalized because they are sufficiently close to their limits, then additional terms can be added to the objective functional in a manner such that the sum of the weights is equal to unity as given in Eq.~\eqref{eq:sum-of-weights}. 
 
\subsection{Modified Optimal Control Problem}\label{sect:modified ocp}

For simplicity, consider the case of a single component of the constraint $c_i,\;i\in[1,\ldots,n_c]$ approaching its limit. Using the penalty term $\C{J}_{P_{\beta}}$ defined in Eq.~\eqref{eq: scaled penalty term}, the optimal control problem defined in Section \ref{sect:ocp} is modified as follows.
Minimize the objective functional 
\begin{equation}\label{eq:modified-ocp-obj}
\C{J}_a = w_1^{(s)} \C{J} + w_2^{(s)} \C{J}_{P_{\beta}}
\end{equation}
subject to the dynamic constraints
\begin{equation}\label{eq:modified-ocp-dyncon}
  \frac{d\m{x}}{dt}-\m{f}(\m{x}(t),\m{u}(t),t)=\m{0}, 
\end{equation}
the inequality path constraints
\begin{equation}\label{eq:modified-ocp-pathcon}
   \m{c}_{\textrm{min}} \leq \m{c}(\m{x}(t),\m{u}(t),t) \leq \m{c}_{\textrm{max}}, \; 
\end{equation}
and the boundary conditions 
\begin{equation}\label{eq:modified-ocp-boundcond}
\begin{array}{rcl}
\m{x}(t_0^{(s)}) &=& \m{x}_0^{(s)}, \\
 \m{b}(\m{x}(t_f^{(s)}),t_f^{(s)}) &=& \m{0}, \\
t_0^{(s)} &= & t_0 + s\Delta T, 
\end{array}
\end{equation}
for $s\in[1,\ldots,S]$. The penalty $\C{J}_{P_{\beta}}$ is not intended to replace a path constraint, but intended to increase the margin between the path constraint function and its limit. Consequently, all path constraints remain in the problem even if a penalty is placed on one or more of the components of the path constraint function.

\subsection{Steps of CCG\&C Method}\label{subsect: steps of constrained cgc}

Suppose that the initial time and initial state are specified, respectively, as $t_0^{(s)}$ and $\m{x}_0^{(s)} $, where $s$ is initialized to zero. The steps of the CCG\&C method are then given as follows:
\begin{enumerate}

\item Set $t_0^{(s)}=t_0+s\Delta T$ and $\m{x}(t_0^{(s)})=\m{x}_0^{(s)}$.
\item Set the weights $(w_1^{(s)},w_2^{(s)})$ as desired for use in Eq.~\eqref{eq:modified-ocp-obj}.
\item Solve the optimal control problem defined by Eqs.~\eqref{eq:modified-ocp-obj}--\eqref{eq:modified-ocp-boundcond} to obtain the solution $\C{B}^{(s)}$.
\item Simulate the perturbed dynamics using the control $\m{u}^*(t)$ obtained in Step (3) over $t\in[t_0^{(s)},t_e^{(s)}]$.
\item Increment $s$ by unity and set $\m{x}_0^{(s)} = \tilde{\m{x}}(t_e^{(s-1)})$.
\item Return to Step (1) until $s=S$.

\end{enumerate}

\section{Application of CCG\&C Method to Reusable Launch Vehicle Entry}\label{Example Application}

This section presents the computationally complex application of reusable launch vehicle (RLV) entry, for which the constrained computational guidance and control method developed in Section~\ref{sect:new-cgs-method} will be implemented. The maximum cross-range RLV entry problem was chosen for this study because the reference solution contains a long-duration constrained arc. The problem formulation in this paper is similar to that of Ref.~\cite{Betts2020}. In this study, state inequality path constraints are placed on the heating rate, dynamic pressure, and sensed acceleration.  This section is organized as follows. First, Section~\ref{subsect:dynamic cons} defines the equations of motion and the inequality constraints that must be satisfied during entry. Next, Section~\ref{subsect:boundary conds} describes the boundary conditions. Finally, Section~\ref{subsect: reference cost} provides the reference objective functional associated with the maximum cross-range problem.

\subsection{Equations of Motion and Path Constraints}\label{subsect:dynamic cons}

\subsubsection{Equations of Motion}\label{subsubsect:equations of motion}
The equations of motion for the RLV entry vehicle (modeled as a point mass over a spherical non-rotating Earth \cite{Betts2020}) are given as
\begin{equation}\label{eq:entry eom}
  \begin{array}{lcl}
    \dot{r} & = & \displaystyle v\sin\gamma, \vspace{3pt} \\
    \dot{\theta} & = & \displaystyle \frac{v\cos\gamma\sin\psi}{r\cos\phi},\vspace{3pt} \\
    \dot{\phi} & =& \displaystyle \frac{v\cos\gamma\cos\psi}{r},\vspace{3pt} \\
    \dot{v} & = & \displaystyle -D-g\sin\gamma , \vspace{3pt} \\
    \dot{\gamma}& =& \displaystyle \frac{1}{v}\left[L\cos\sigma+\left(\frac{v}{r}-\frac{g}{v}\right)\cos\gamma\right] , \vspace{3pt} \\
    \dot{\psi}& =& \displaystyle \frac{L\sin\sigma}{v\cos\gamma}+\frac{v\cos\gamma\sin\psi\tan\theta}{r}, 
  \end{array}
\end{equation}
where $r$ is the geocentric radius, $v$ is the speed, $\theta$ is the longitude, $\phi$ is the latitude, $\gamma$ is the flight path angle,  $\psi$ is the azimuth (measured clockwise from due north), $\sigma$ is the bank angle, $g = \mu/r^2$ is the gravitational acceleration, and $\mu$ is Earth's gravitational parameter. The lift and drag specific forces are defined as
\begin{equation}\label{eq:lift and drag specific forces}
  \begin{array}{lcl}
    D & =& \displaystyle qSC_D/m , \vspace{3pt} \\
    L & = & \displaystyle qSC_L/m, 
  \end{array}
\end{equation}
where $q = \rho v^2/2$ is the dynamic pressure, $\rho=\rho_0\exp(-h/H)$ is the atmospheric density,  $h = r - R_e$ is the altitude, $R_e$ is the radius of the Earth, $\rho_0$ is the surface level atmospheric density, $H$ is the density scale height, $C_D$ is the coefficient of drag, $C_L$ is the coefficient of lift, $S$ is the vehicle reference area, and $m$ is the vehicle mass. The coefficients of drag and lift are modeled as
\begin{equation}\label{eq: coeff lift and drag}
  \begin{array}{lcl}
    C_D & =& \displaystyle C_{L_0}+C_{L_1} \alpha , \vspace{3pt} \\
    C_L & = & \displaystyle C_{D_0}+C_{D_1} \alpha + C_{D_2} \alpha^2, 
  \end{array}
\end{equation}
where $\alpha$ is the angle-of-attack and $(C_{D_0},C_{D_1}, C_{D_2}, C_{L_0}, C_{L_1})$ are constants. Finally, in this study, the angle of attack and bank angle rates, denoted $u_\alpha$ and $u_\sigma$, respectively, are used as controls.  Consequently, the following two differential equations are augmented to the system given in Eq.~\eqref{eq:entry eom}:
\begin{equation}\label{eq: control rates}
  \begin{array}{lcl}
    \dot{\alpha} & = & \displaystyle u_{\alpha}, \vspace{3pt} \\
    \dot{\sigma} & =& \displaystyle u_{\sigma} .
  \end{array}
\end{equation}

Including $\alpha$ and $\sigma$ as states ensures  continuity in the angle of attack and bank angle across guidance cycles as well as providing more realistic variations in these quantities throughout the motion. The state is then given as $\m{x}(t)=\left[r(t),\theta(t),\phi(t),v(t),\gamma(t),\psi(t),\alpha(t),\sigma(t)\right]'$, and the control as $\m{u}(t) = \left[u_\alpha(t), u_\sigma(t)\right]'$.

\begin{table}[ht]
\centering
\caption{Nominal physical and aerodynamic constants \label{tab:physical-constants-Earth}}
\renewcommand{\baselinestretch}{1.25}\normalsize\normalfont
\begin{tabular}{ccc}\hline\hline
Quantity & Value & Unit  \\\hline
  $R_e$ & $6.371\times10^6$ & m   \\
  $r_n$ & $1$ & m \\
  $\mu$ & $3.986\times10^{14}$ & $\textrm{m}^3$/$\textrm{s}^2$ \\
  $g_0$ & $9.80665$ & m/$s^2$ \\
  $S$ & $249.9$ & $\textrm{m}^2$ \\
   $K_q$ & $1.7415\times10^{-4}$ & W/$\textrm{cm}^2$ \\
  $C_{D_0}$ & $0.0785$ & - \\
   $C_{D_1}$ & $-6.1593\times 10^{-3}$ & 1/deg \\
    $C_{D_2}$ & $6.2142\times10^{-4}$ & 1/deg$^2$  \\
  $C_{L_0}$ & $-0.207$ & - \\
  $C_{L_1}$ & $0.029245$ & 1/deg \\
  $m$ & $92079$ & kg  \\
  $\rho_0$ & $1.2256$ & kg/$\textrm{m}^3$\\
   $H$ & $7254$ &m  \\\hline\hline
\end{tabular}
\end{table}

\subsubsection{Inequality Constraints}\label{subsect: inequality cons}

The following constraints are placed on the angle of attack and bank angle: 
\begin{equation}\label{eq: control bounds}
  \begin{array}{lclcl}
    \alpha_{\min} & \leq &  \displaystyle \alpha & \leq & \alpha_{\max} , \vspace{3pt} \\
    \sigma_{\min} & \leq& \displaystyle \sigma & \leq & \sigma_{\max}  \vspace{3pt}.
  \end{array}
\end{equation}
Next, the angle of attack and bank angle rates are bounded by 
\begin{equation}\label{eq: control rates bounds}
  \begin{array}{lclcl}
    u_{\alpha,{\min}} & \leq &  \displaystyle u_\alpha & \leq & u_{\alpha,{\max}} , \vspace{3pt} \\
    u_{\sigma,{\min}} & \leq &  \displaystyle u_\sigma & \leq & u_{\sigma,\max} .\vspace{3pt} 
  \end{array}
\end{equation}
Additionally, three state-inequality path constraints are enforced during entry. First, the stagnation point heating rate \cite{sutton1971general} is constrained by
  \begin{equation}\label{eq: heating rate con}
  \begin{array}{lclcl}
    \dot{Q}_{\min} & \leq& \displaystyle \dot{Q} & \leq & \dot{Q}_{\max}  \vspace{3pt},
  \end{array}
\end{equation}
where  
\begin{equation}\label{eq:heating rate}
 \dot{Q} = K_q \sqrt{\frac{\rho}{r_n}}v^3,
\end{equation}
$K_q$ is a heating rate constant, and $r_n$ is the nose radius of the entry vehicle. Second, the following constraint is enforced on the dynamic pressure
 \begin{equation}\label{eq: dynamic pressure}
  \begin{array}{lclcl}
    q_{\min} & \leq& \displaystyle q & \leq & q_{\max}  \vspace{3pt}.
  \end{array}
\end{equation}
Third, the sensed acceleration is constrained by
  \begin{equation}\label{eq: sensed accel con}
  \begin{array}{lclcl}
    A_{\min} & \leq& \displaystyle A & \leq & A_{\max}  \vspace{3pt},
  \end{array}
\end{equation}
where 
 \begin{equation}\label{eq:sensed accel}
      A = \frac{1}{g_0}\sqrt{L^2+D^2}  
 \end{equation}
 is the sensed acceleration normalized by $g_0$, the standard surface level gravitational acceleration of Earth. Together, the path constraint function given by $\m{c}(t) = [\dot{Q}(t),q(t),A(t)]'$ will be the focus for the constraint violation prevention method described in Section~\ref{sect:new-cgs-method}. The nominal aerodynamic coefficients and physical parameters used in this study are given in Table \ref{tab:physical-constants-Earth}.

\begin{table}[ht]
\centering
\caption{Constraint limits \label{tab:constraint limits}}
\renewcommand{\baselinestretch}{1.25}\normalsize\normalfont
\begin{tabular}{cccc}\hline\hline
  Quantity & Unit & Minimum & Maximum\\\hline
  $\alpha$ & deg &$-90$  & $90$ \\
  $\sigma$ & deg & $-90$  & $1$ \\
  $u_\alpha$ & deg/s &  $-0.5$  & $0.5$ \\
   $u_\sigma$ & deg/s & $-5$  & $5$ \\
  $\dot{Q}$ & MW/m$^2$ & $0$  & $0.85$ \\
  $q$ & kPa & $0$  & $16.375$\\
  $A$ & $ g $ & $0$  & $2.5$\\\hline\hline
\end{tabular}
\end{table}
 
\subsection{Boundary Conditions}\label{subsect:boundary conds}
The initial and terminal conditions are given by
\begin{equation}\label{eq:boundary condition}
\begin{array}{lclclcl}
    r(t_0^{(s)}) & = & r_0^{(s)} & , & r(t_f^{(s)}) & = & r_f,    \\
    \theta(t_0^{(s)}) & = & \theta_0^{(s)} & , & \theta(t_f^{(s)}) & = & \textrm{Free}, \\
    \phi(t_0^{(s)}) & = & \phi_0^{(s)}  & , & \phi(t_f^{(s)}) & = & \textrm{Free}, \\
    v(t_0^{(s)}) & = & v_0^{(s)}  & , & v(t_f^{(s)}) & = & v_f, \\
    \gamma(t_0^{(s)}) & = & \gamma_0^{(s)} & , & \gamma(t_f^{(s)}) & = & \gamma_f, \\
    \psi(t_0^{(s)}) & = & \psi_0^{(s)} & , & \psi(t_f^{(s)}) & = & \textrm{Free}, \\
    \alpha(t_0^{(s)}) & = & \alpha_0^{(s)} & , & \alpha(t_f^{(s)}) & = & \textrm{Free}, \\
    \sigma(t_0^{(s)}) & = & \sigma_0^{(s)} & , & \sigma(t_f^{(s)}) & = & \textrm{Free}.
    \end{array}
\end{equation}
At the start of each guidance cycle $s$, the initial conditions are set to the values of the actual state at the guidance update time, as described by Section~\ref{subsect: steps of constrained cgc}. The terminal conditions, however, are held constant at each guidance update. The initial conditions at the specified initial time, $t_0$, (corresponding to $s = 0$), and terminal conditions are given in Table \ref{tab:bcs-Earth}.

\begin{table}[ht]
\centering
\caption{Boundary conditions \label{tab:bcs-Earth}}
\renewcommand{\baselinestretch}{1.25}\normalsize\normalfont
\begin{tabular}{cccc}\hline\hline
  Quantity & Unit & Initial & Terminal\\\hline
  $r$ & km &$645$  & $639.6$ \\
  $\theta$ & deg & $0$  & Free \\
  $\phi$ & deg&  $0$  & Free \\
   $v$ & m/s & $7802.88$  & $762$ \\
  $\gamma$ & deg & $-1$  & $-5$ \\
  $\psi$ & deg & $90$  & Free\\ \hline\hline
\end{tabular}
\end{table}

\subsection{Reference Objective Functional}\label{subsect: reference cost}
The reference objective functional, $\C{J}$, of the RLV entry problem is as follows. Maximize the final cross-range, which is equivalently formulated as minimizing the negative of the final latitude 
\begin{equation}\label{eqn:rlv obj}
\C{J} = -\phi(t_f).
\end{equation}
The reference optimal control problem can now be stated as follows: minimize the objective functional given by Eq.~\eqref{eqn:rlv obj}, subject to the dynamic constraints of Section~\ref{subsubsect:equations of motion}, the inequality constraints of Section~\ref{subsect: inequality cons}, and the boundary conditions of Section~\ref{subsect:boundary conds}. In Section~\ref{subsec:comparing weights}, the objective functional given by Eq.~\eqref{eqn:rlv obj} is modified as described by Section~\ref{sect:new-cgs-method}, for the purpose of preventing path constraint violations. Then, in Section~\ref{subsec:Reducing phugoids} an additional term is incorporated into the objective functional, intended to reduce phugoid oscillations.

\section{Results of CCG\&C Method on RLV Entry Problem}\label{results}

In this section, the constrained computational guidance and control method developed in Section~\ref{sect:new-cgs-method} is applied to the RLV entry problem formulated in Section~\ref{Example Application}, and its performance is compared against the computational guidance and control method of Ref.~\cite{Dennis2019}. This section is organized as follows. First, in Section~\ref{Nominal Reference Solution} the key features of the reference solution of the RLV entry problem formulated in Section~\ref{Example Application} are provided. Second, the CCG\&C method is employed for a single value of perturbed surface level density and the modified objective functional formulated in Section~\ref{sect:new-cgs-method} is studied by varying the objective functional weights. Third, in Section~\ref{subsec:Reducing phugoids} an additional penalty term, which is intended to reduce phugoid oscillations, is included in the modified objective functional. Finally, in Section~\ref{subsect:Monte-Carlo simulations} a Monte-Carlo campaign consisting of 1000 perturbed cases is used to verify and compare the performance of the CCG\&C method, with and without the aforementioned phugoid penalty, against the performance of the CG\&C method of Ref.~\cite{Dennis2019}.

The results shown in Sections~\ref{Nominal Reference Solution}-\ref{subsect:Monte-Carlo simulations} were obtained via the following software/hardware. All optimal solutions were generated using the MATLAB optimal control software $\mathbb{GPOPS-II}$ \cite{Garg2010,Garg2011April,Garg2011June,Patterson2015}, which uses LGR collocation to transcribe the continuous time optimal control problem into a large sparse NLP. $\mathbb{GPOPS-II}$ was employed with the NLP solver IPOPT \cite{Biegler2008} in full Newton mode, with the initial offline solution obtained using relative error tolerance set to $10^{-7}$ and maximum iterations set to 2000. First and second derivatives were supplied using sparse central finite differencing. Mesh refinement was employed using the method of Ref.~\cite{Liu2018}, with a mesh refinement error tolerance of $10^{-5}$, the maximum number of mesh refinement iterations set to 10, and the allowable collocation points per interval set to a minimum of three and maximum of 20. When re-solving the optimal control problem to make guidance updates, the mesh refinement error tolerance was increased to $10^{-4}$, and a maximum of two mesh refinement iterations was allowed. Once the optimal control problem is solved at the start of a guidance cycle, a cubic spline interpolant is used to provide a continuous approximation for the control over the simulated horizon using the MATLAB ordinary differential equation solver $\textsf{ode113}$ with a relative error tolerance of $10^{-8}$. All computations were performed on a 2023 Apple M3 MacBook Pro running macOS Sonoma 14.3.1 with 96 GB RAM. Simulations were run in MATLAB version R2023b.

\subsection{Reference Solution }\label{Nominal Reference Solution}

\begin{figure}[ht]
\centering

\subfloat[ Heating rate, $\dot{Q}(t)$ vs time, $t$ \label{fig:heating rate nominal}]{\includegraphics[scale=0.4]{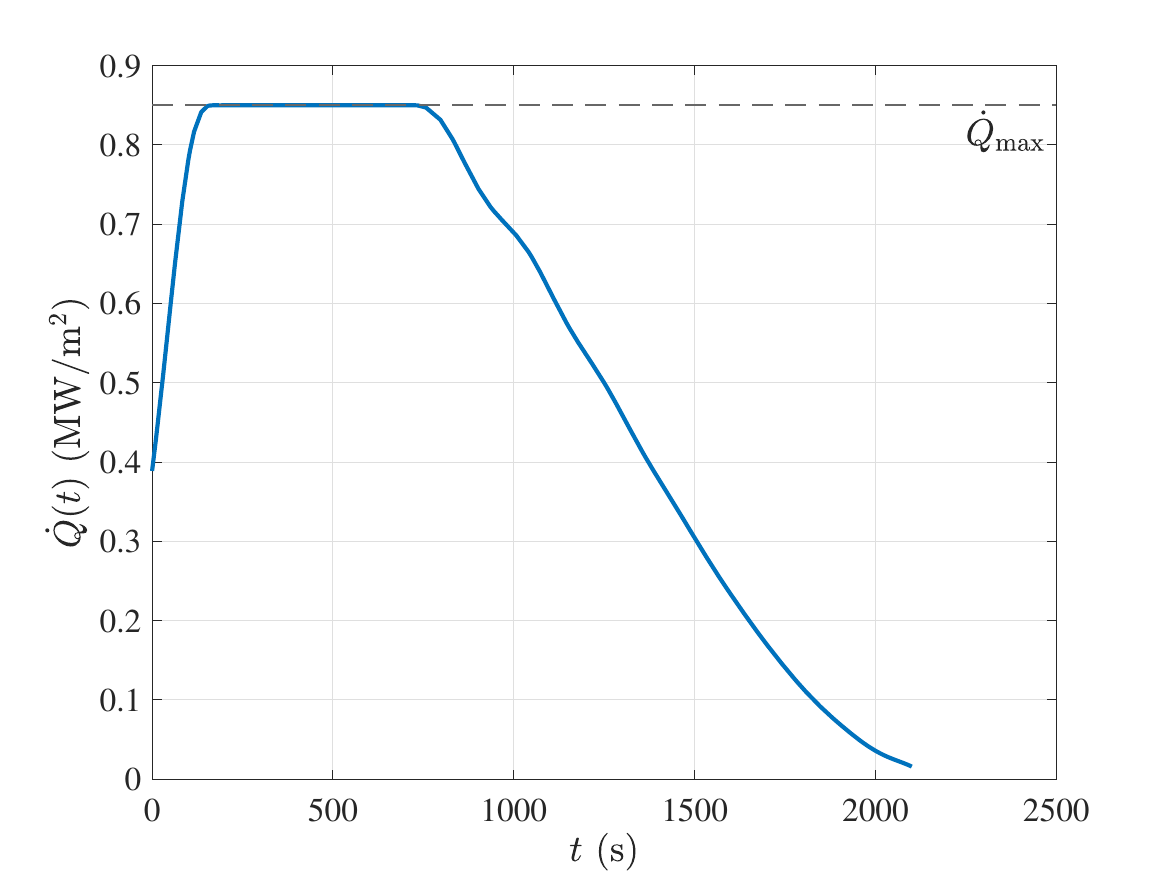}}~~~\subfloat[ Dynamic Pressure, $q (t) $ vs time, $t$\label{fig:dynamic pressure nominal}]{\includegraphics[scale=0.4]{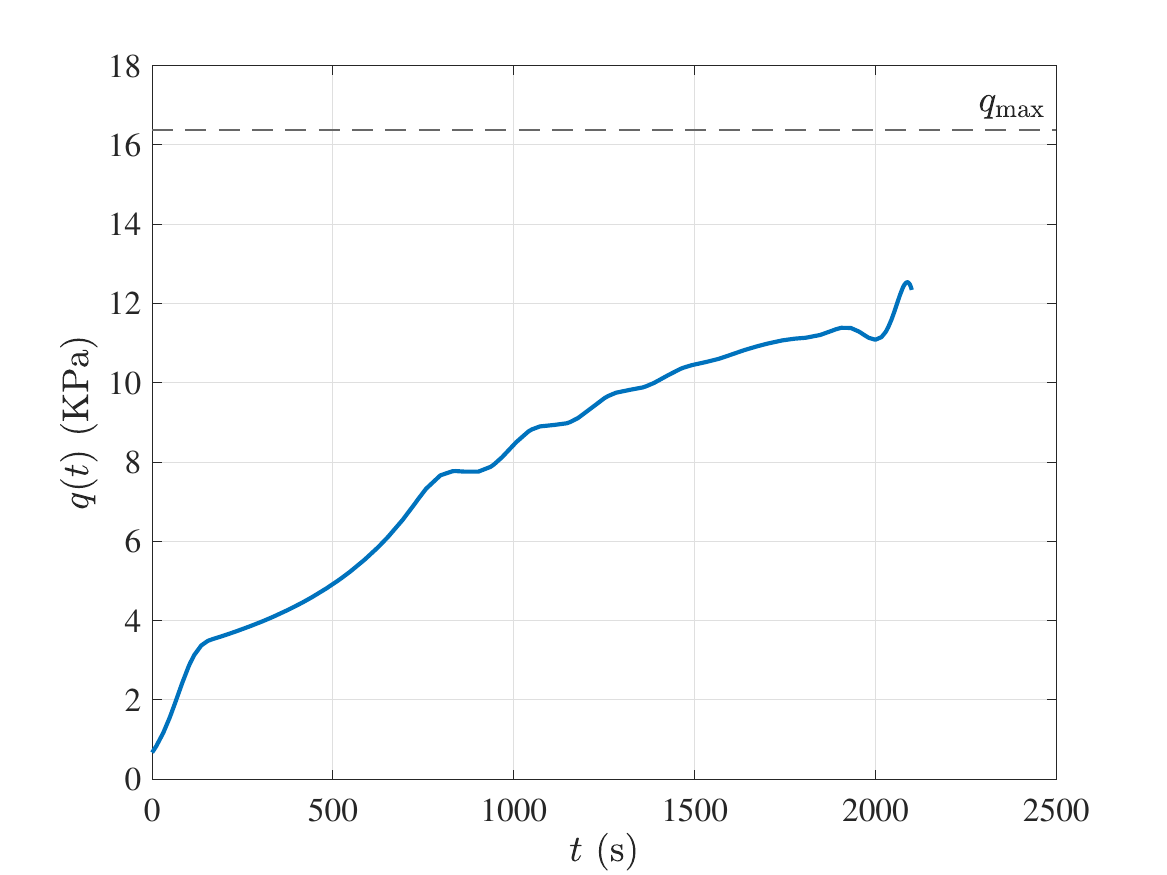}}

\subfloat[ Sensed acceleration, $ A (t) $ vs time, $t$ \label{fig:sensed accel nominal}]{\includegraphics[scale=0.4]{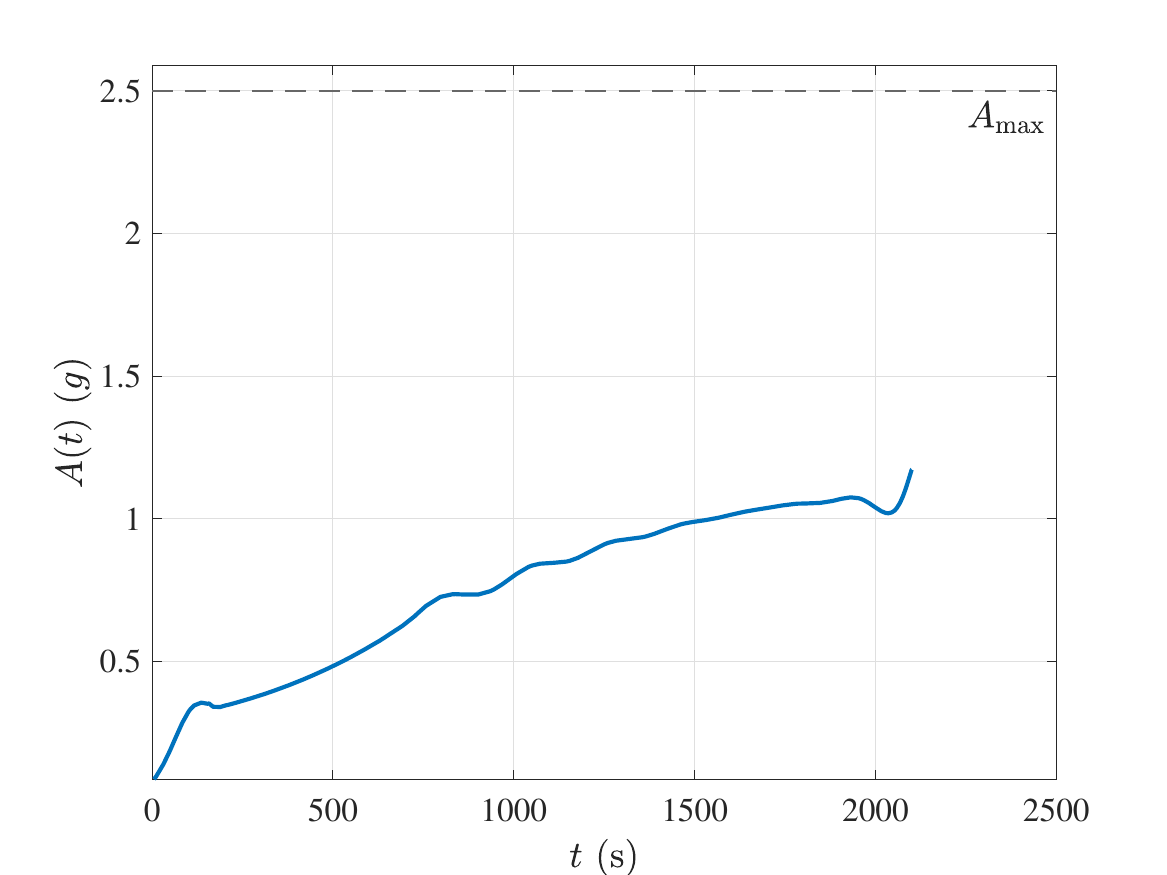}}

  \caption{Nominal path constraint profiles.}\label{fig:path cons nominal}
\end{figure}

Key characteristics of the reference solution of the maximum cross range RLV entry problem presented in Section~\ref{Example Application} are studied here. First, the optimal objective on the reference solution is computed as $\phi^*(t_f) \approx 34.0$~deg. Second, from the constraint profiles shown in Fig.~\ref{fig:path cons nominal}, it is seen that the reference solution contains a long-duration heating rate constrained arc segment. The constrained arc occurs on the interval $t \in [168.69,725.87]$, accounting for approximately $26.5\%$ of the entire trajectory duration.  Third, after the vehicle exits the constrained arc, the heating rate decreases with time for the remainder of the trajectory. Without any compensation, deviations from the reference optimal solution may violate the path constraint and the optimal control problem will become infeasible. It is noted that the dynamic pressure and sensed acceleration constraint functions maintain a large margin between their maximum values and their limits and do not become active on the reference solution.  Nevertheless, these constraints remain in the problem formulation to ensure the NLP is feasible.

\clearpage

\subsection{Application of CCG\&C method}\label{subsec:comparing weights}

In this section, the CCG\&C method developed in Section~\ref{sect:new-cgs-method} is applied to the RLV entry problem using a perturbed surface level atmospheric density $\tilde{\rho}_0 = 1.01 \rho_0$. Furthermore, the CCG\&C method is employed using the following parameter values: $\beta=5$ as required in Eq.~\eqref{log-integral-exponential-scaled} and  $\eta_i = c_{\max,i}$ as required in Eq.~\eqref{eq: scaled penalty term}.  The performance of the CCG\&C method for these values is then analyzed in terms of the objective functional weights, and compared against the CG\&C method of Ref.~\cite{Dennis2019}. For both methods, guidance updates are made every $\Delta T = 15$~s. When the remaining horizon of the simulated trajectory is too small to obtain a viable solution to the NLP, the optimal control computed from the previous guidance update is used for the remainder of the simulation. 

Recall from Section~\ref{Nominal Reference Solution} that, with the exception of the segment where the heating rate constraint is active, the path constraint functions are not in close proximity of their limits.  

As a result, it is necessary to add a penalty term into the objective functional only when a path constraint lies within a user-specified threshold of a path constraint limit at the terminus of a guidance cycle.  
 
The decision to add this penalty can then be made using the following threshold function

\begin{equation}\label{eq: weights on guid}
    (w_1^{(s)},w_2^{(s)}) = \left\{\begin{array}{lcl} (1,0) & , &    \frac{{c}_i\left(t_e^{(s-1)}\right)}{c_{\max,i}}  < \xi , \\ (\zeta_1,\zeta_2) & , &   \frac{{c}_i\left(t_e^{(s-1)}\right)}{c_{\max,i}}  \geq \xi, \end{array} \right. 
\end{equation}
where $\xi$ is a user-specified threshold parameter, and the first two arguments of $c_i(\tilde{\m{x}}(t_e^{(s-1)}),\m{u}(t_e^{(s-1)}),t_e^{(s-1)})$ are omitted for brevity. As only the heating rate constraint attains its limit, it is the only path constraint function for which the threshold function of Eq.~\eqref{eq: weights on guid} is relevant. 

For the results shown here, a threshold of $\xi = 0.9$ was employed. The performance of the CCG\&C method is now studied for varying values of $(\zeta_1,\zeta_2)$. Note that the results shown for $(\zeta_1,\zeta_2) = (1,0)$ correspond to the CG\&C method of Ref.~\cite{Dennis2019}, for which the optimal control problem is re-solved using only the reference objective functional.

\begin{figure}[h]
\centering

\subfloat[ Altitude, $h$ vs speed, $v$ \label{fig:alt versus speed weight comps}]{\includegraphics[scale=0.4]{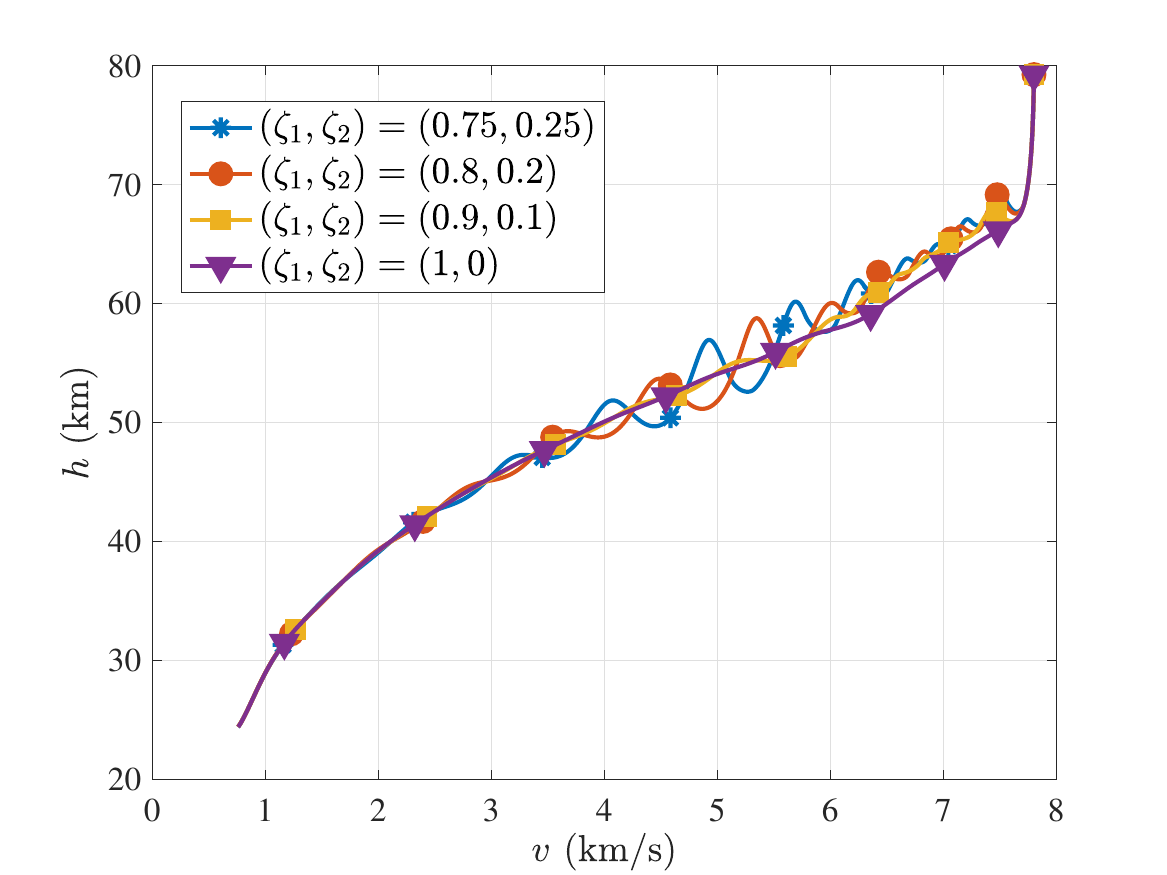}}~~~\subfloat[ Latitude, $ \phi  $ vs longitude, $\theta$ \label{lat v lon weigths compared}]{\includegraphics[scale=0.4]{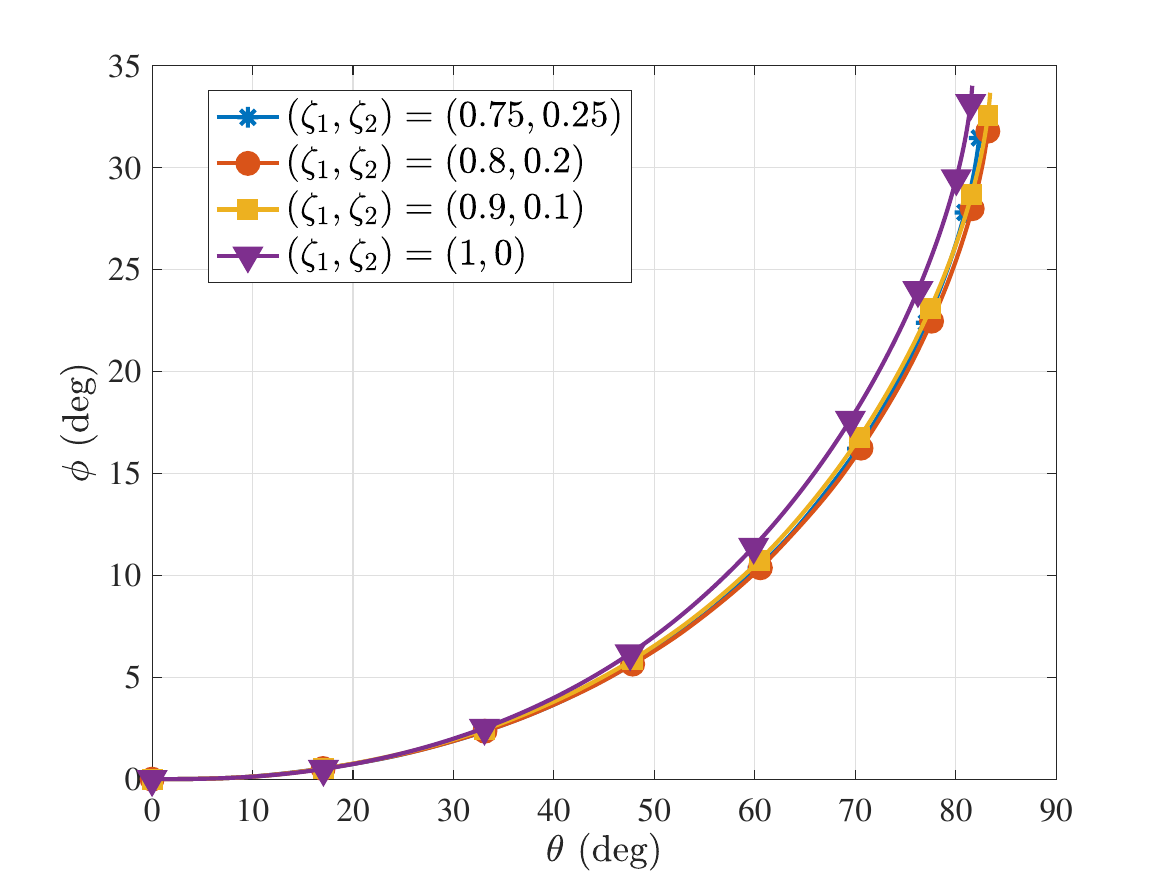}}

\subfloat[ Flight path angle, $ \gamma (t) $ vs time, $t$ \label{fpa weigths compared}]{\includegraphics[scale=0.4]{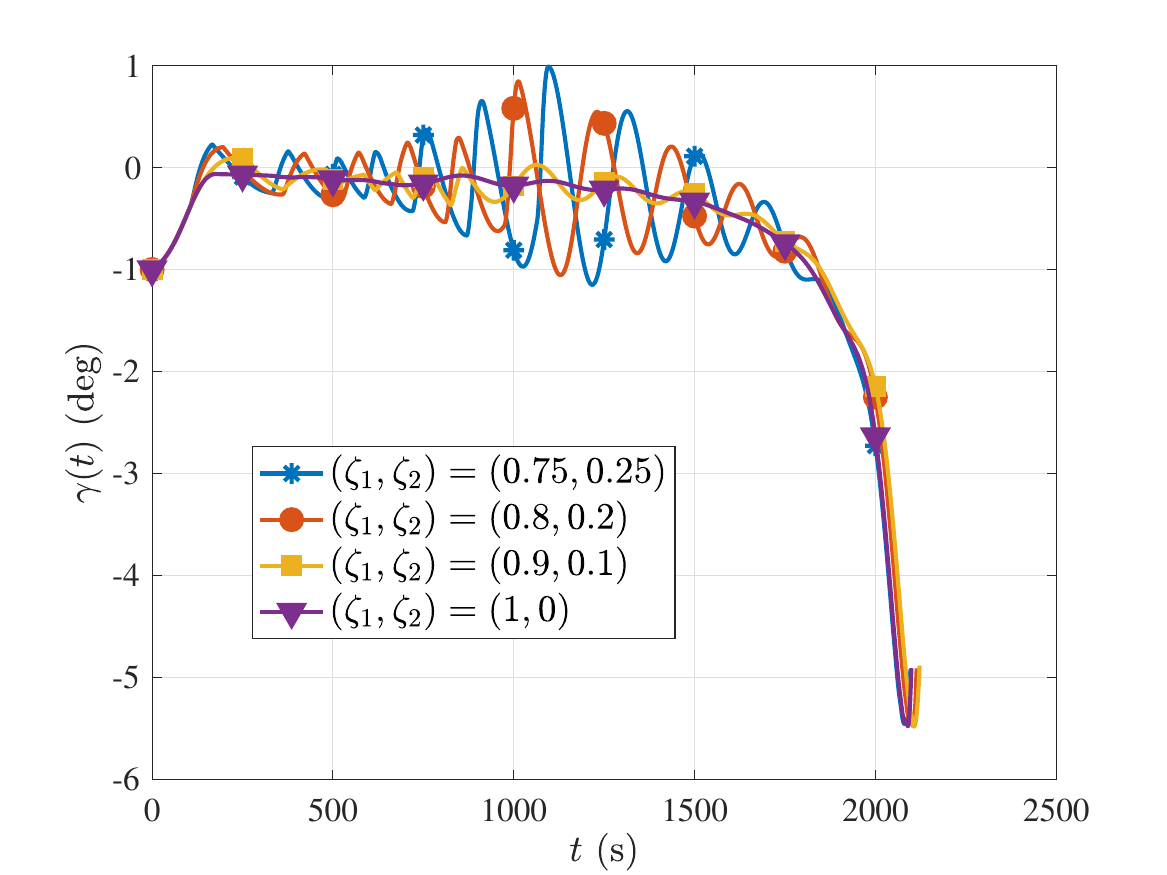}}

  \caption{State profiles obtained from simulation of perturbed dynamics for different values of $(\zeta_1,\zeta_2)$. }\label{fig: states weights compared}
\end{figure}

\begin{figure}[h]
\centering

\subfloat[ Heating rate, $\dot{Q}(t)$ vs time, $t$ \label{fig:heating rate weights compared full view}]{\includegraphics[scale=0.4]{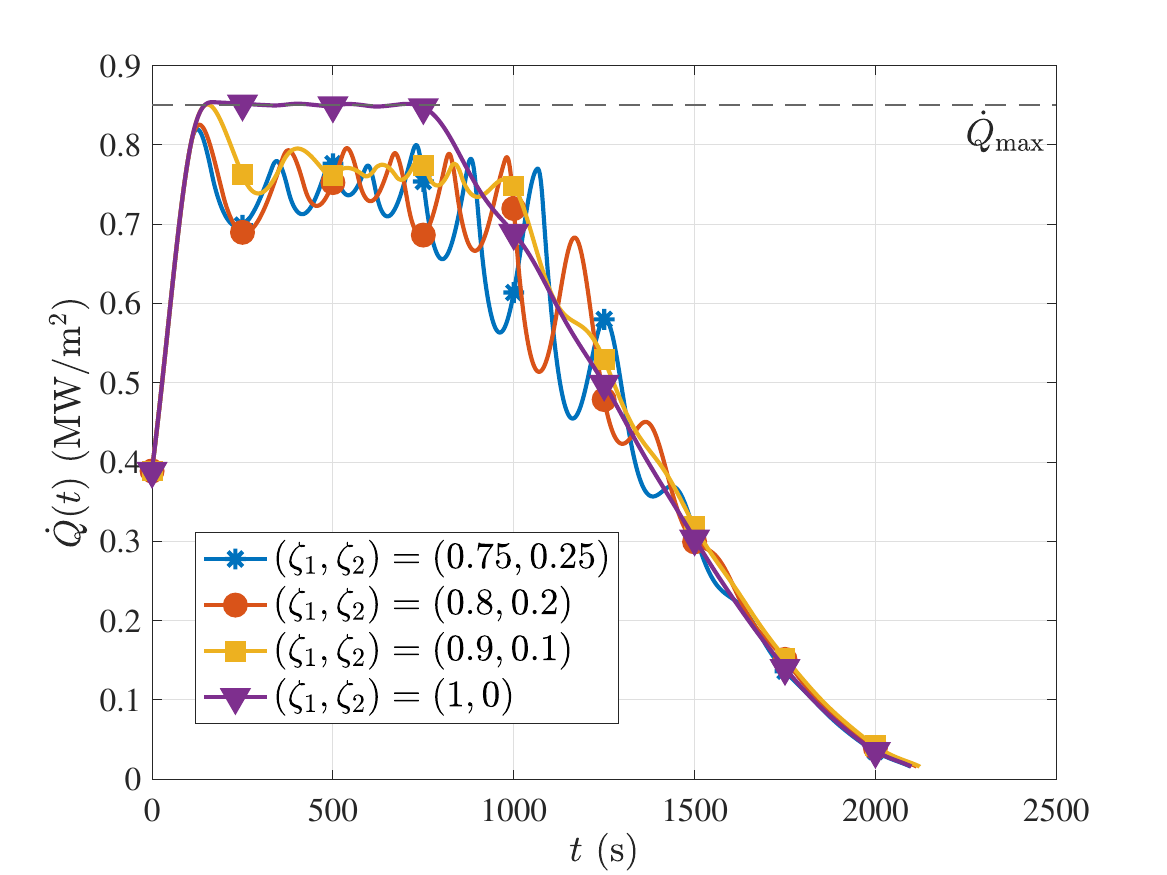}}~~~\subfloat[ Enlarged view of constraint violations\label{fig:heating rate enlarged}]{\includegraphics[scale=0.4]{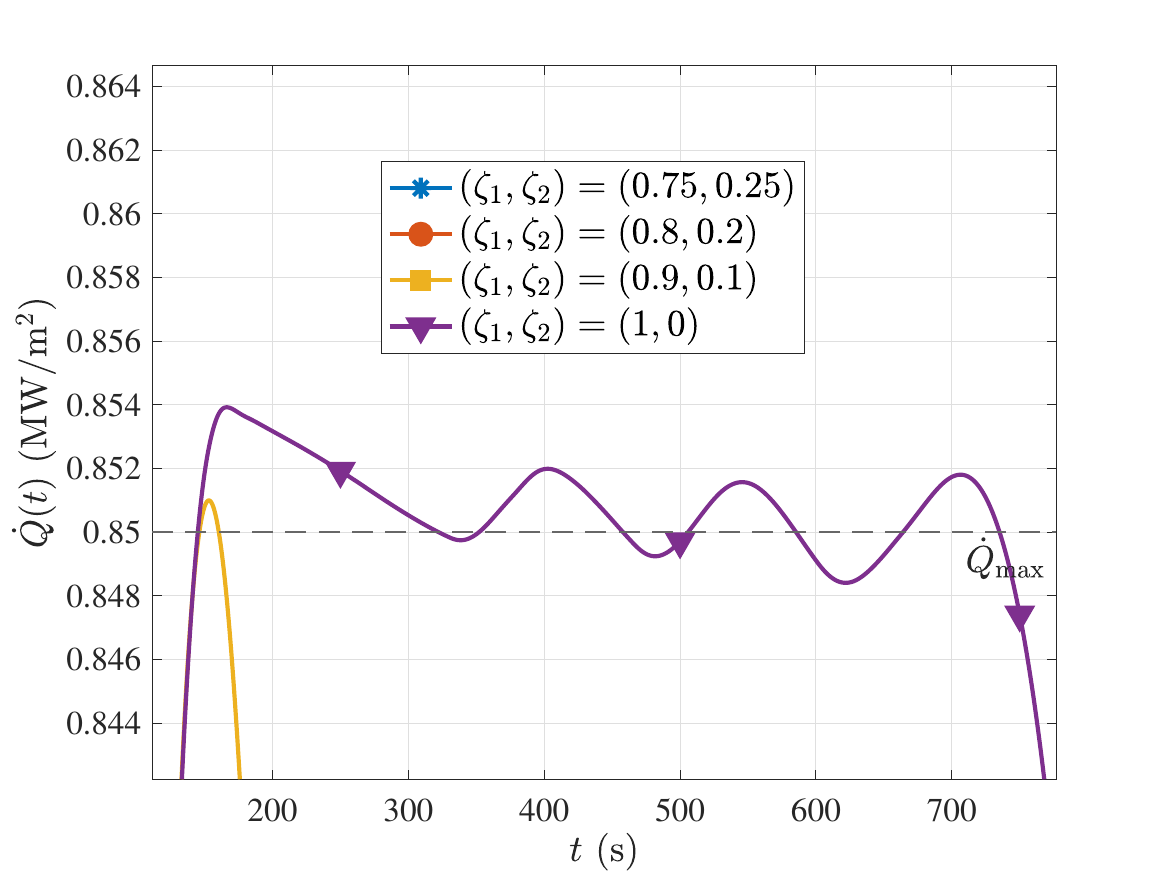}}

  \caption{Heating rate profiles obtained from simulation of perturbed dynamics for different values of $(\zeta_1,\zeta_2)$. }\label{fig:heating rate weights compared}
\end{figure}

\begin{figure}[h]
\centering

\subfloat[ Dynamic pressure, $q (t) $ vs time, $t$\label{fig:dynamic pressure weights compared}]{\includegraphics[scale=0.4]{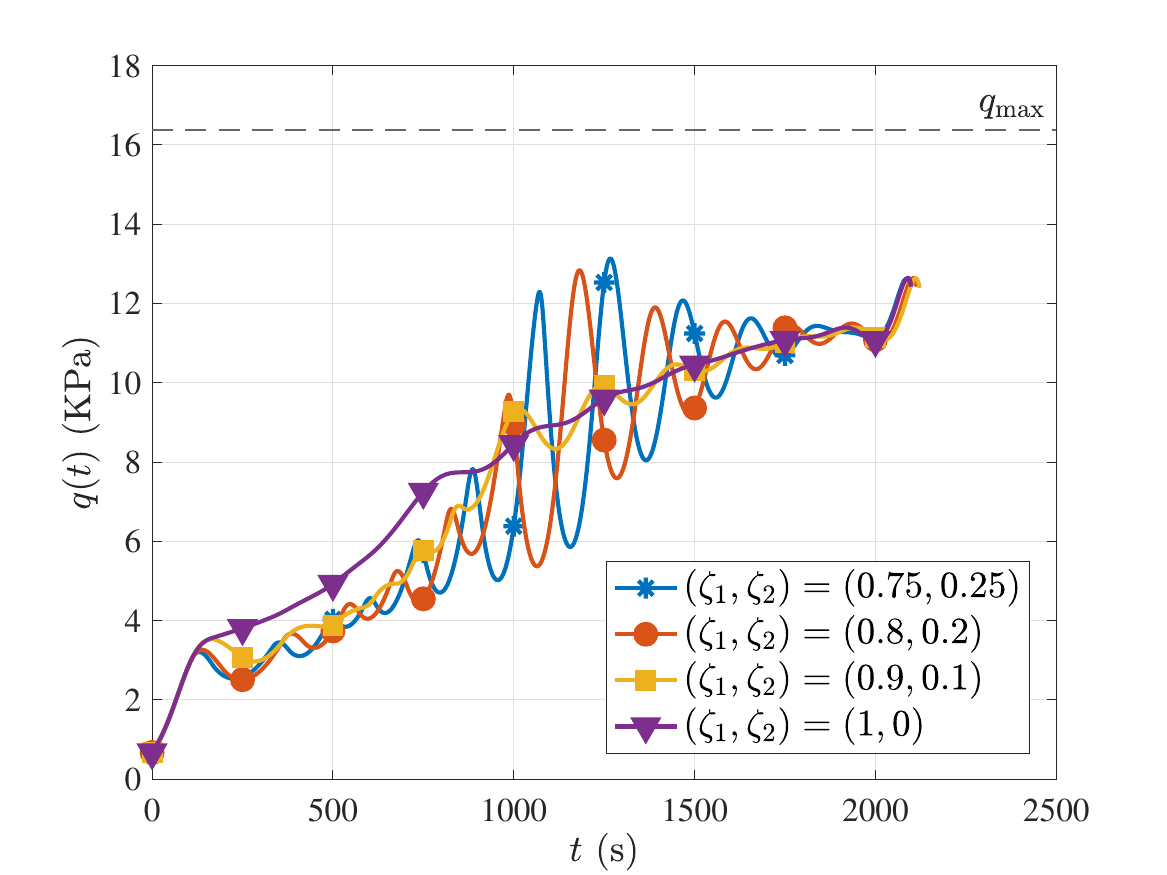}}~~~\subfloat[ Sensed acceleration, $ A (t) $ vs time, $t$ \label{fig:sensed accel weights compared}]{\includegraphics[scale=0.4]{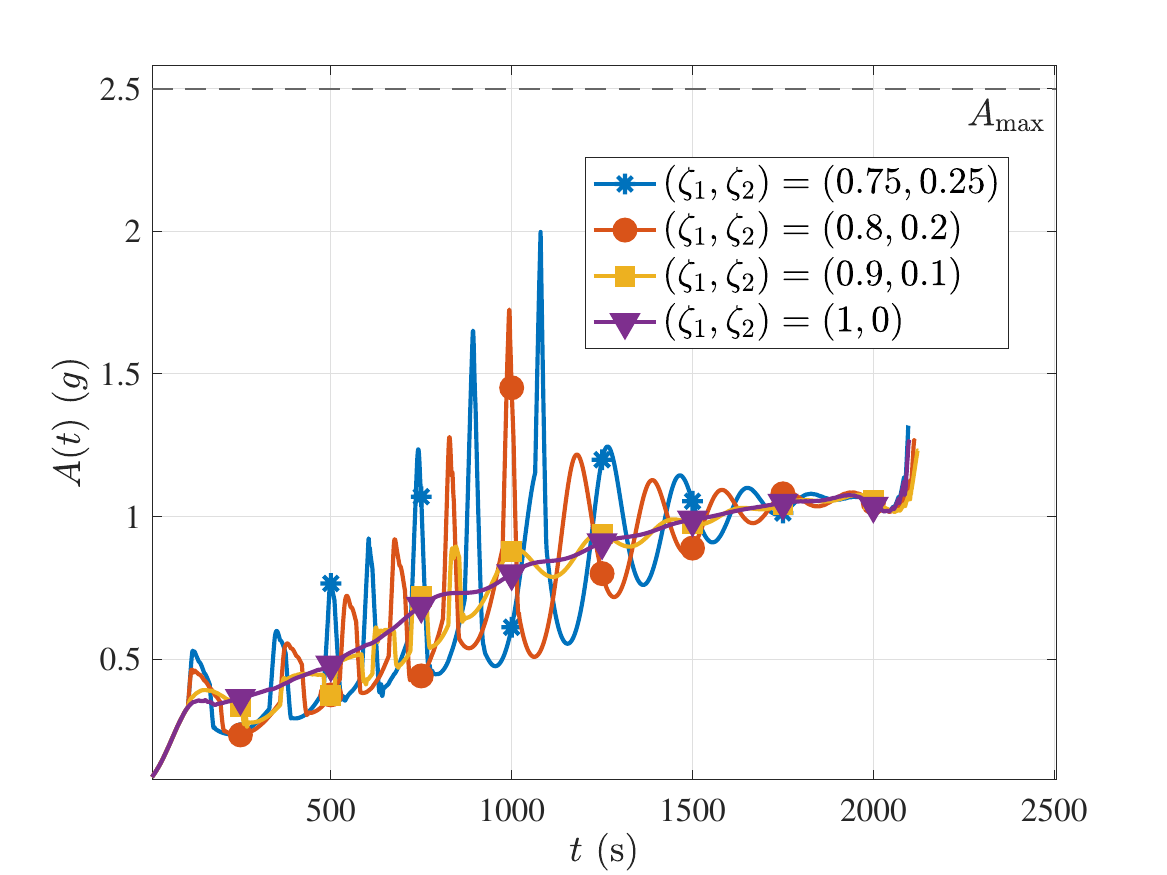}}

  \caption{Dynamic pressure and sensed acceleration profiles obtained from simulating perturbed dynamics for different values of $(\zeta_1, \zeta_2)$. }\label{fig:path cons weights compared}
\end{figure}

\begin{figure}[h]
\centering

\subfloat[ Angle of attack, $\alpha(t)$ vs time, $t$ \label{fig:aoa weights compared}]{\includegraphics[scale=0.4]{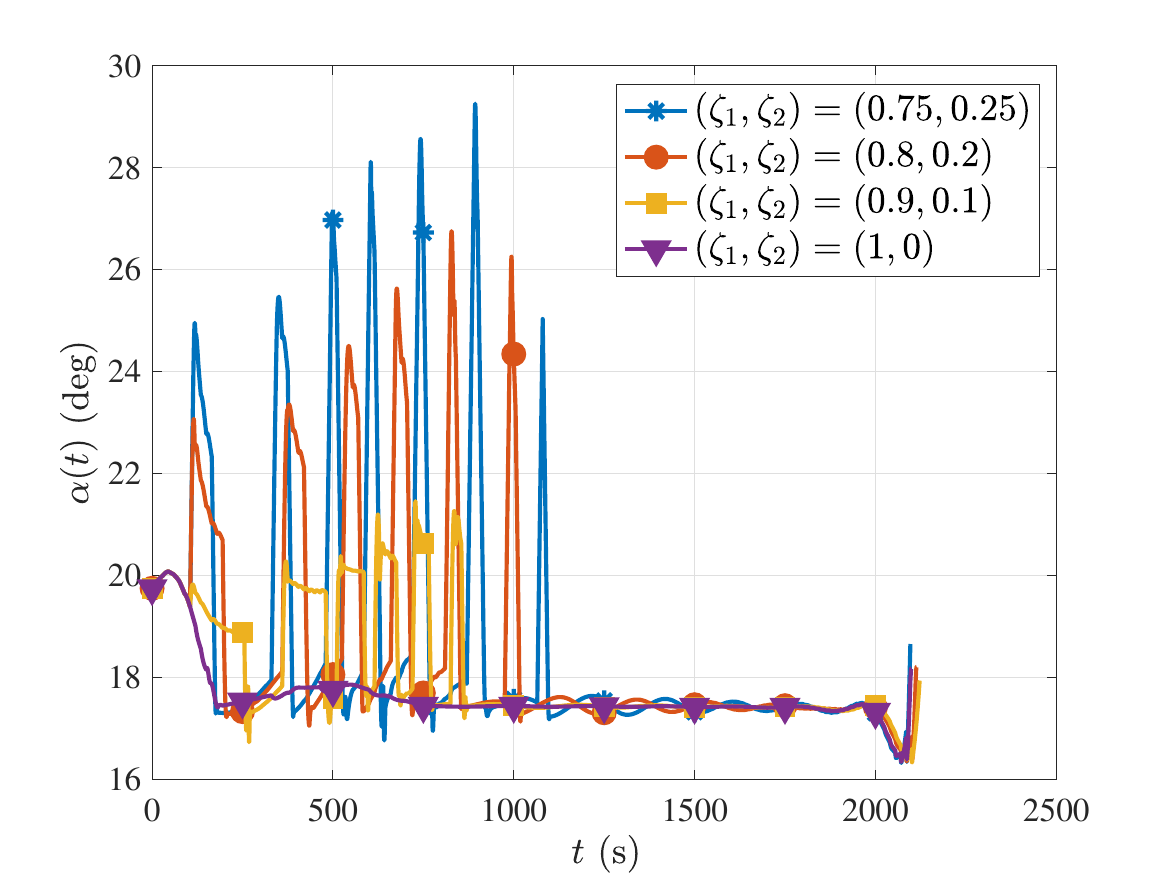}}~~~\subfloat[ Bank angle, $\sigma (t) $ vs time, $t$\label{fig:sigma weights compared}]{\includegraphics[scale=0.4]{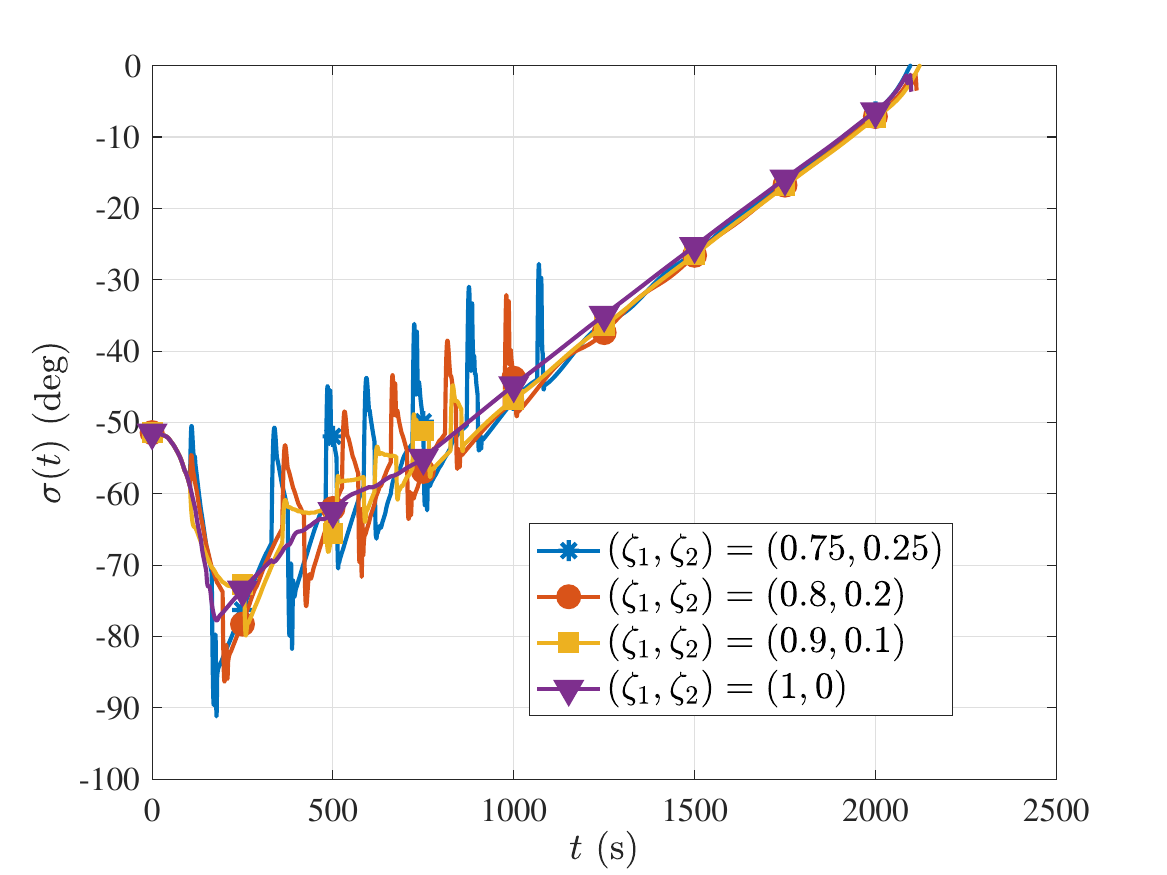}}

  \caption{Angle of attack and bank angle profiles for different values $(\zeta_1,\zeta_2)$. }\label{fig:controls weights comp}
\end{figure}

Figure \ref{fig:heating rate enlarged} shows the heating rate for four different combinations of $(\zeta_1,\zeta_2)$.  It is seen that, for the case $(\zeta_1,\zeta_2)=(1,0)$, the heating rate profile of the perturbed system contains several segments where the constraint limit is violated. During guidance cycles where the heating rate constraint is violated, the optimal control problem is infeasible, meaning guidance updates cannot be performed. Note that when the optimal control problem is infeasible, the control obtained on the last feasible solution is used until an update can be obtained. As a result of constraint violations, a total of $30$ guidance cycles were infeasible using the CG\&C method of Ref.~\cite{Dennis2019}.  Now, consider the trajectories obtained with the CCG\&C method. Observing the perturbed heating rate profiles in Fig.~\ref{fig:heating rate weights compared}, it can be seen that when the path constraint penalty is included, the constraint violations are significantly reduced. In particular, the duration of constraint violations decrease as the weight on the penalty term in the objective function increases.  For the case of $(\zeta_1,\zeta_2) = (0.9,0.1)$, where the weight on the penalty term is small, a brief constraint violation starting at $t \approx 146$~s occurs, resulting in a single guidance cycle where an update cannot be performed. As $\zeta_2$ is increased, however, the maximum value of the heating rate is reduced such that the path constraint is never violated.

Table \ref{tab: weights comp results} summarizes the key numerical results obtained from varying the values of $(\zeta_1,\zeta_2)$ used in the modified objective functional. A trade-off exists between maximizing the final latitude $\phi (t_f)$ and lowering the peak value of the heating rate. When employing the CG\&C method, a final latitude is reached that is in agreement with the reference solution, but this result comes at the cost of path constraint violations. On the other hand, as the constraint function penalty weight $\zeta_2$ is increased, the final latitude decreases. While some performance is lost when penalizing the path constraint, simultaneously feasibility is improved.  In particular, the margin between the maximum heating rate of the perturbed system increases as $\zeta_2$ increases. It is also noted that the errors in the terminal value of the state are all on the same order of magnitude and lie within reasonable ranges.

While penalizing the path constraint either reduces or eliminates infeasible guidance cycles, Fig.~\ref{fig:path cons weights compared} shows that, large amplitude oscillations are seen in the sensed acceleration and dynamic pressure profiles.  In particular,
Fig.~\ref{fig:controls weights comp} demonstrates that the angle of attack and bank angle change abruptly for any guidance cycle where even a small penalty on the path constraint function is included, and these large changes in angle of attack and bank angle induce phugoid oscillations. Figure \ref{fig: states weights compared} shows that once the phugoid oscillations start, they remain in the solution for almost the entire trajectory (except towards the very end of the trajectories).  Such phugoid oscillations are known to be undesirable during atmospheric entry \cite{miller2022end,vinh1981phugoid,liu2016rapid}. In the next section, an additional penalty is included in the objective function to assist in attenuating phugoid oscillations.

\begin{table}[h]
\centering
\caption{Numerical results obtained for varying weights $(\zeta_1,\zeta_2)$ used in modified objective functional. \label{tab: weights comp results}}
\renewcommand{\baselinestretch}{1.25}\normalsize\normalfont
\begin{tabular}{c c c c c c }\hline\hline
  $(\zeta_1,\zeta_2)$ &$\phi(t_f)$  & $\delta h (t_f)$  & $\delta v (t_f)$  & $\delta \gamma (t_f)$ & $\max (\dot{Q})$ \\
   & (deg) & (m) & (m/s) & (deg) & (MW/m$^2$) \\\hline
  $(0.75,0.25)$ & $32.270$ &$6.351$ & $-0.6287$ &  $0.0757$ & $0.8198$ \\
   $(0.8,0.2)$ & $32.785$ &$9.010$ & $-0.7479$ & $0.0870$ & $0.8254$ \\
   $(0.9,0.1)$ & $33.659$ & $20.341$ & $-1.105$ & $0.1144$ & $0.8510$ \\
   $(1,0)$ & $34.004$ &$9.412$ & $-0.7642$ & $0.08846$ & $0.8539$ \\\hline\hline
\end{tabular}
\end{table}

\subsection{Reducing Phugoid Oscillations}\label{subsec:Reducing phugoids}

In Section~\ref{subsec:comparing weights}, the constraint function penalty in the modified objective functional, given by Eq.~\eqref{eq:modified-ocp-obj}, reduced path constraint violations. These solutions, however, were prone to phugoid oscillations.  This section builds upon the previous section by augmenting the modified objective functional given by Eq.~\eqref{eq:modified-ocp-obj} with an additional penalty term on the flight path angle, where this additional term is intended to reduce phugoid oscillations. The phugoid penalty used in this study employs an approach similar to that given in Ref.~\cite{miller2022end}.

The results shown here were obtained using the following modified objective functional
\begin{eqnarray}\label{modified obj 2}
\C{J}_a = w_1^{(s)} \C{J} + w_2^{(s)}\C{J}_{P_{\beta}} +  w_3^{(s)} \int_{t_0}^{t_f} \left( \frac{C}{1+e^{-k \sin \gamma}}-\frac{C}{2} \right)^2 dt,
\end{eqnarray}
where $C$ is a constant defined as
\begin{equation}
    C = 2 \frac{1+e^{-k}}{1-e^{-k}}
\end{equation}
such that the value of the integrand in Eq. \eqref{modified obj 2} is equal to zero at $\gamma = 0$ deg, and monotonically increases to a value of unity when approaching values of $\gamma = \pm 90$ deg. The user chosen constant $k$ dictates how quickly the value within the integrand approaches unity as $\gamma$ strays from zero degrees; as $k$ increases, the phugoid penalty grows faster as $\gamma$ deviates from 0 degrees. This study is performed using $k=3$, which was chosen as this value appropriately penalizes deviations in the flight path angle without overly restricting the motion of the vehicle. The weights $(w_1^{(s)},w_2^{(s)},w_3^{(s)})$ were set at the start of each guidance cycle according to the threshold function
\begin{equation}\label{eq: phu weights on guid}
    (w_1^{(s)},w_2^{(s)},w_3^{(s)}) = \left\{\begin{array}{lcl} (1,0,0) & , &   \frac{{c}_i\left(t_e^{(s-1)}\right)}{c_{\max,i}}  < \xi , \\ (\zeta_1,\zeta_2,\zeta_3) & , &   \frac{{c}_i\left(t_e^{(s-1)}\right)}{c_{\max,i}}  \geq \xi.\end{array} \right. 
\end{equation}
 The phugoid penalty is only included in the objective functional when the path constraint penalty was included, so as to prevent the phugoid oscillations induced by the path constraint penalty as shown in Section~\ref{subsec:comparing weights}. It was desired to find a value of $\zeta_3$ large enough to prevent phugoid oscillations, but not so large as to impede the path constraint penalty from sufficiently increasing the margin between the constraint function and its limit. In this section, the results are shown for the perturbed surface level density $\tilde{\rho}_0 = 1.01\rho_0$, and were obtained using the parameters $(\zeta_1,\zeta_2) = (0.8,0.2)$, $\beta = 5$, $\eta_i = c_{\max,i}$, $\xi = 0.9$. The effect of the additional phugoid penalty term is studied for values of $\zeta_3= (1,3,5)$.

 \begin{figure}[h]
\centering

\subfloat[ Altitude, $h$ vs speed, $v$ \label{fig:alt versus speed wphu comps}]{\includegraphics[scale=0.4]{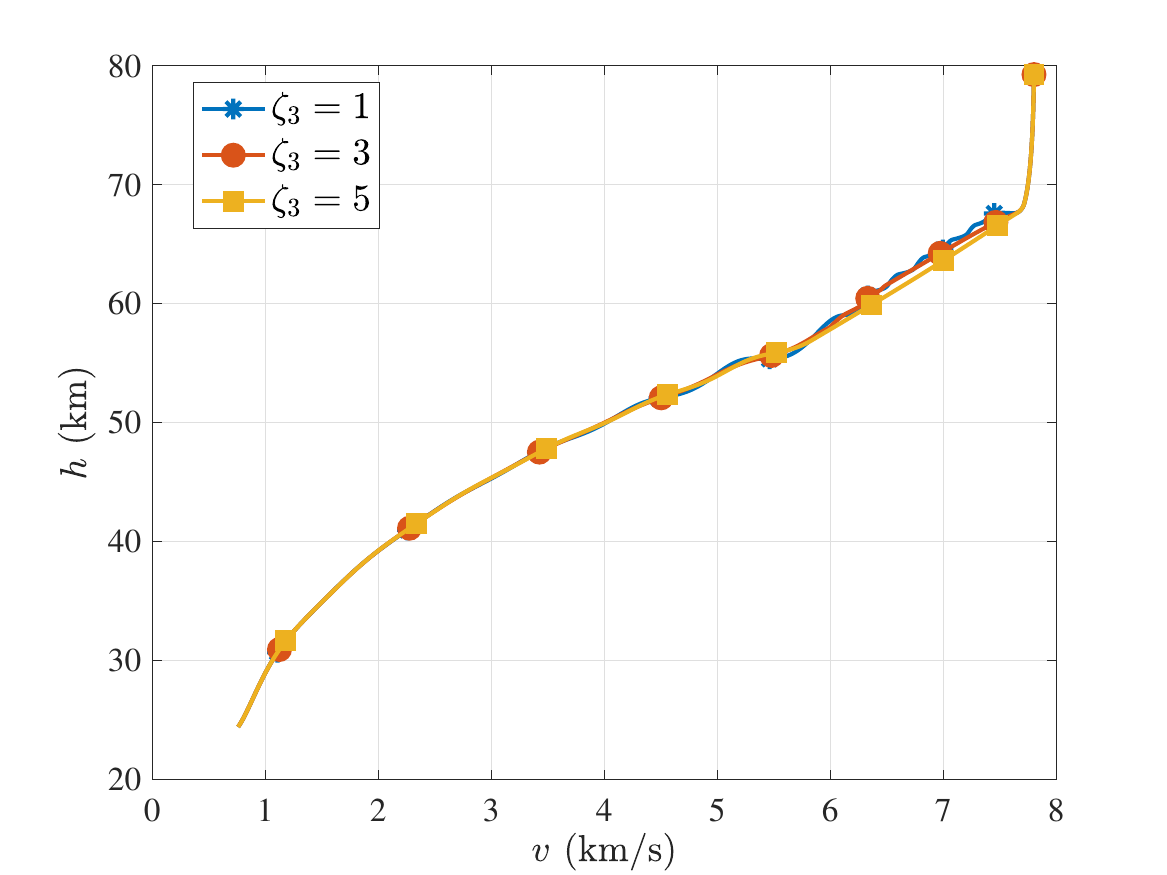}}~~~\subfloat[ Latitude, $ \phi $ vs longitude, $\theta$ \label{lat v lon phu compared}]{\includegraphics[scale=0.4]{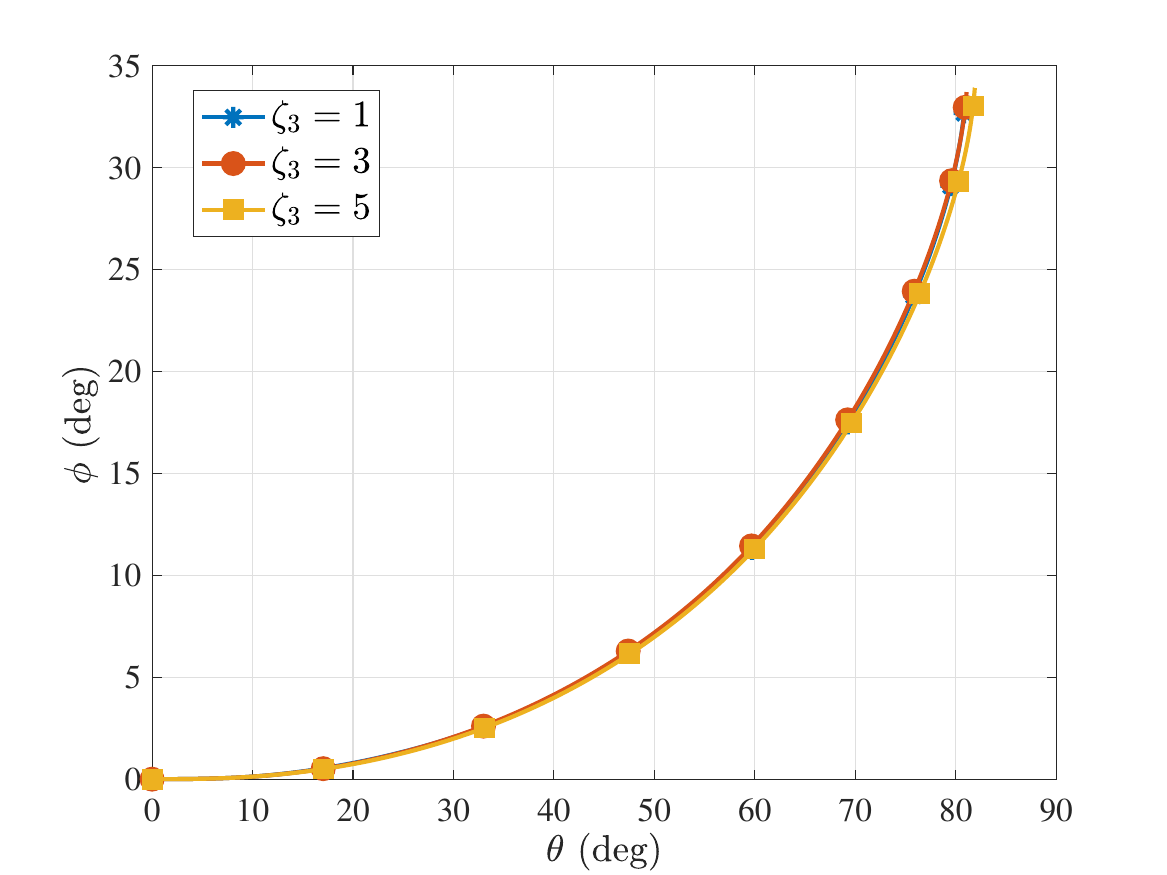}}

\subfloat[ Flight path angle, $ \gamma (t) $ vs time, $t$ \label{fpa phu compared}]{\includegraphics[scale=0.4]{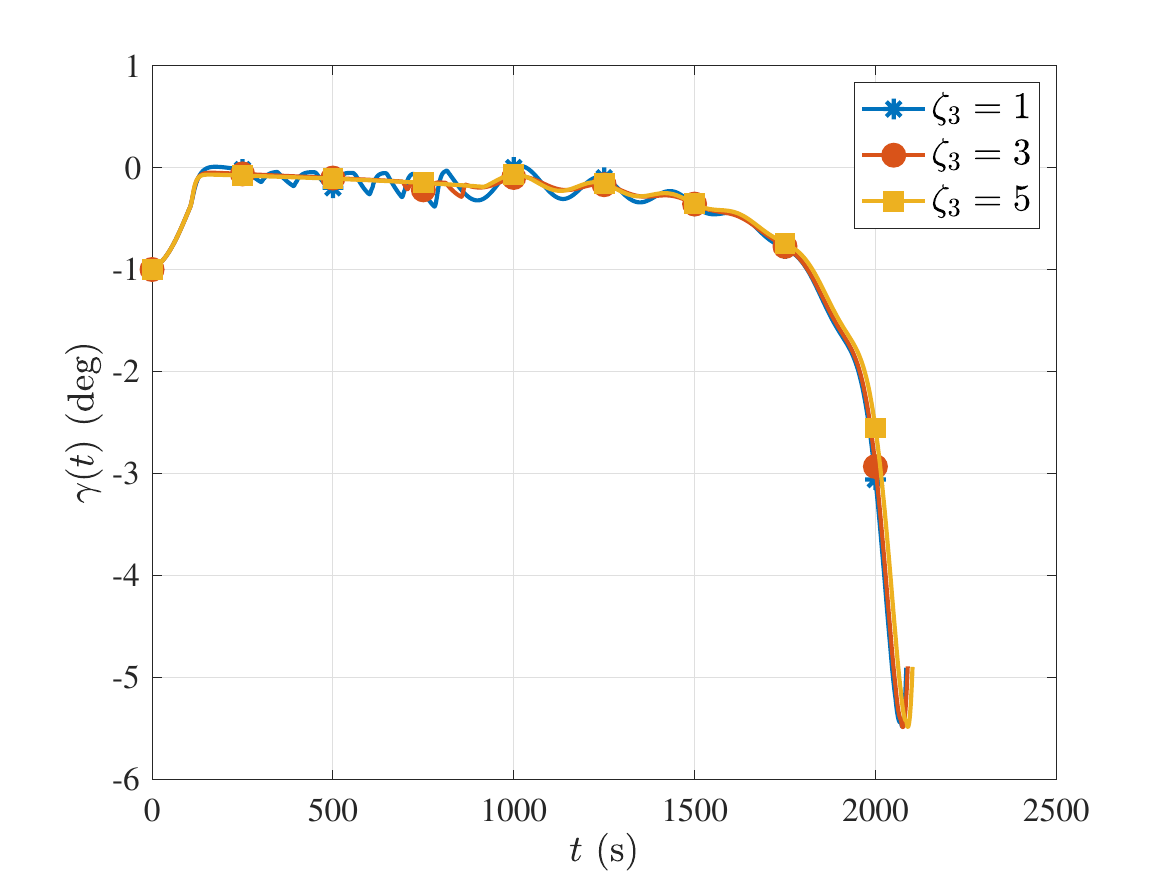}}

  \caption{State profiles obtained from simulation of perturbed dynamics for varied values of $\zeta_3$. }\label{fig: states phu compared}
\end{figure}

\begin{figure}[h]
\centering

\subfloat[ Heating rate, $\dot{Q}(t)$ vs time, $t$ \label{fig:heating rate phu compared}]{\includegraphics[scale=0.4]{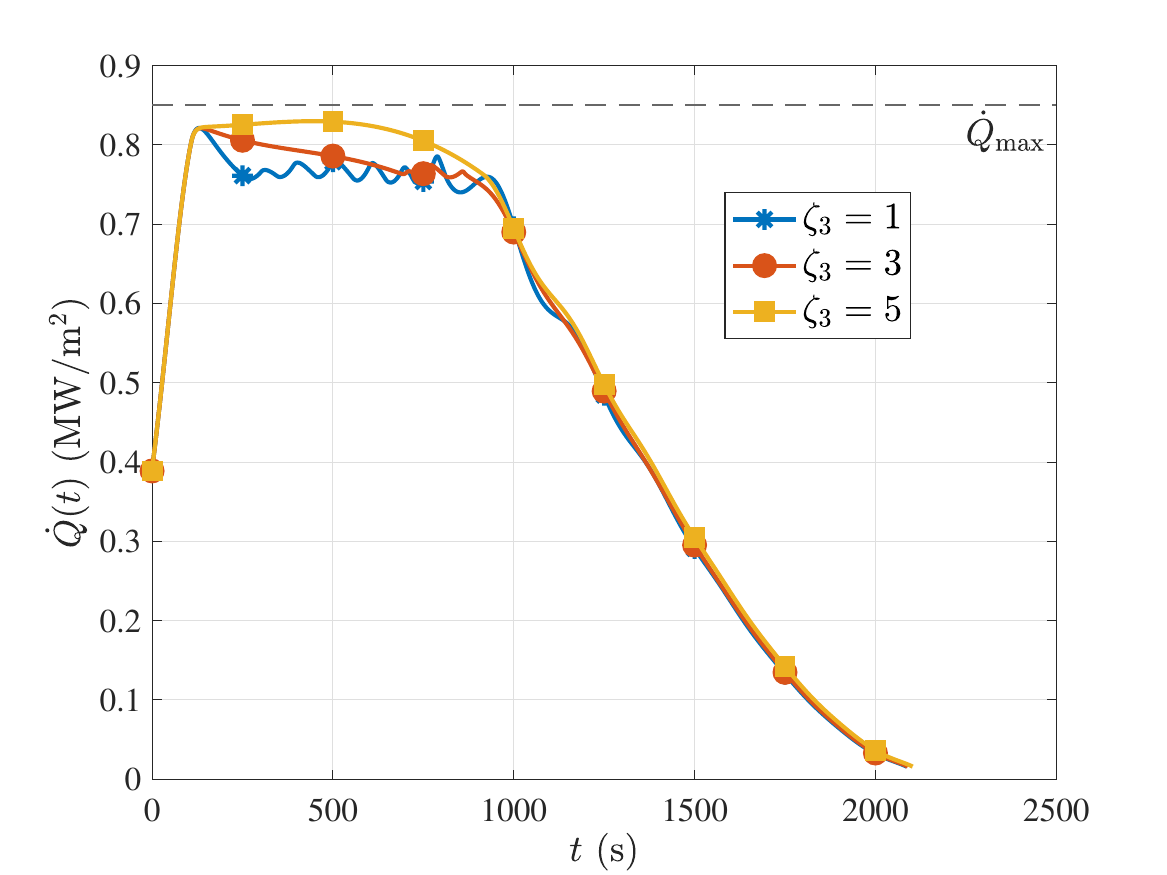}}~~~\subfloat[ Dynamic pressure, $q (t) $ vs time, $t$\label{fig:dynamic pressure phu compared}]{\includegraphics[scale=0.4]{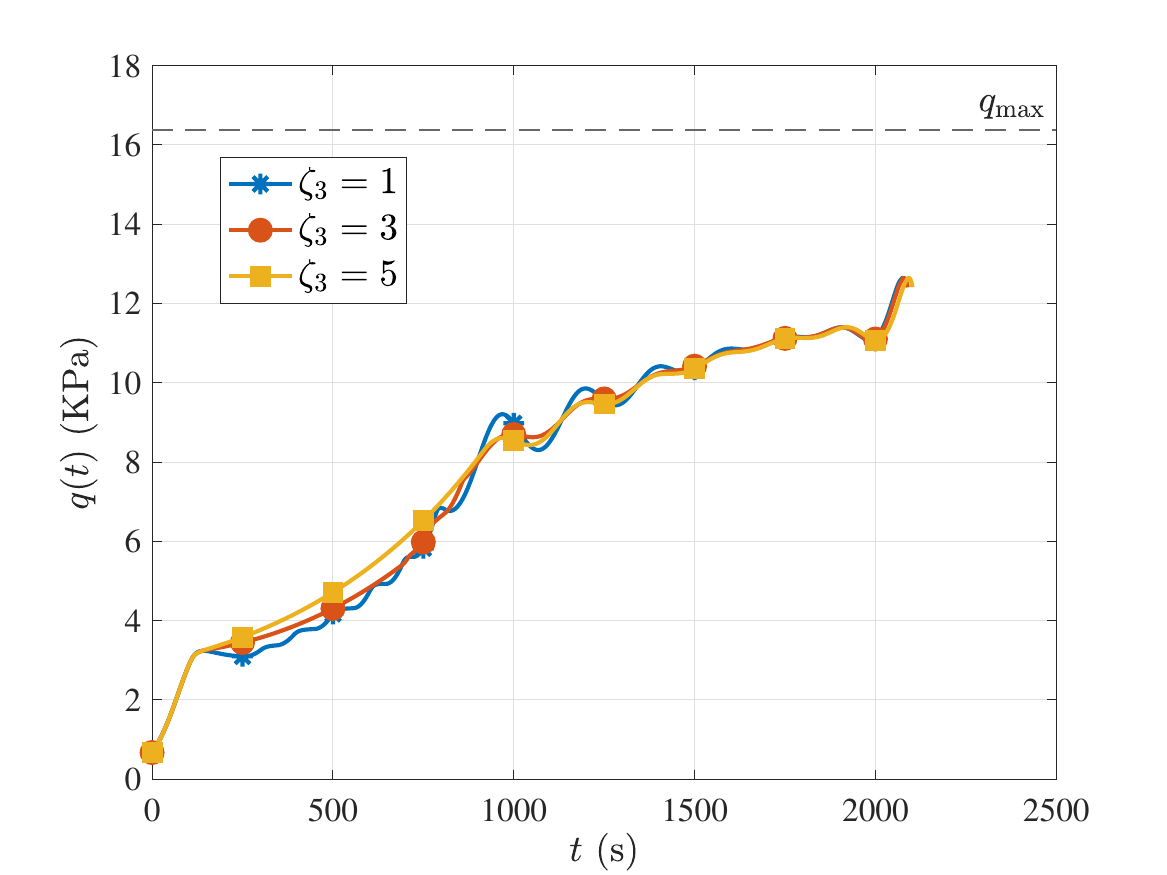}}

\subfloat[ Sensed acceleration, $ A (t) $ vs time, $t$ \label{fig:sensed accel phu compared}]{\includegraphics[scale=0.4]{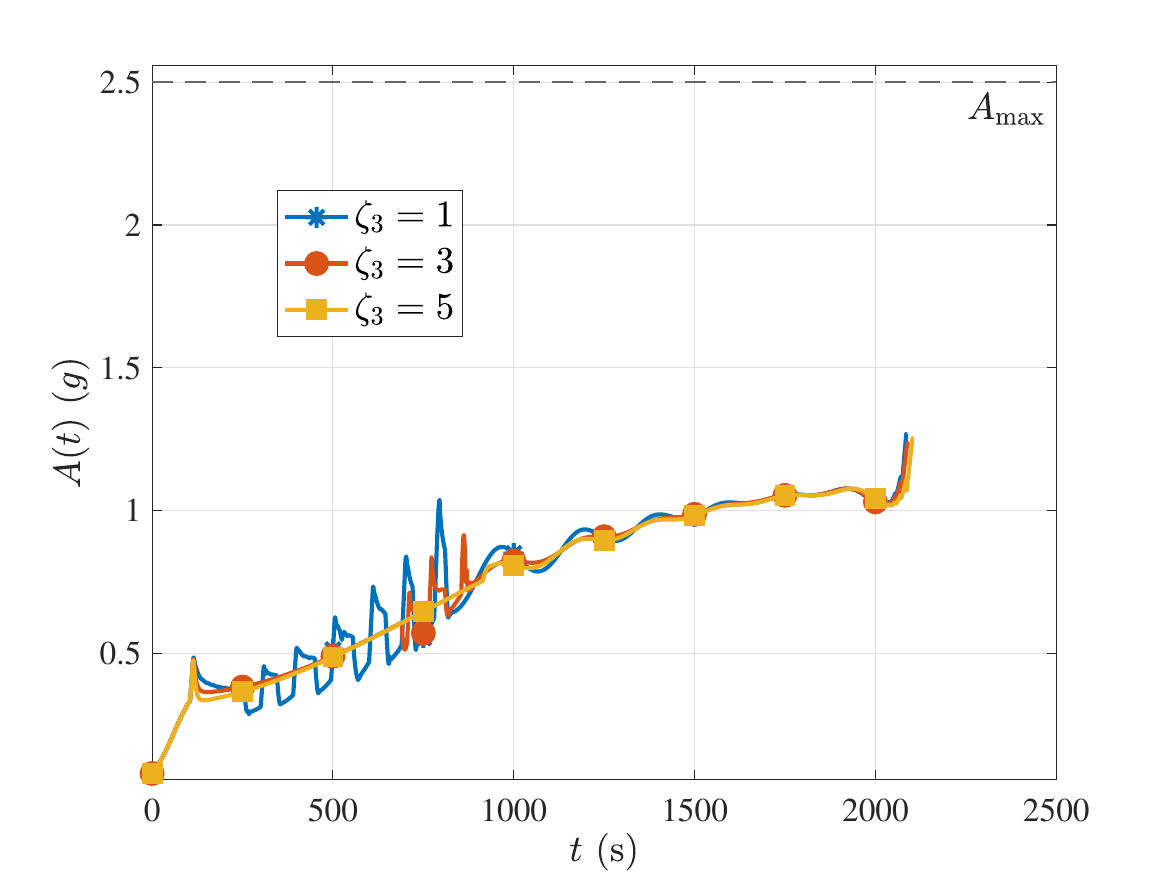}}

  \caption{Path constraint profiles obtained from simulation of perturbed dynamics for varied values of $\zeta_3$. }\label{fig:path cons phu compared}
\end{figure}

\begin{figure}[h]
\centering

\subfloat[ Angle of attack, $\alpha(t)$ vs time, $t$ \label{fig:aoa phu compared}]{\includegraphics[scale=0.4]{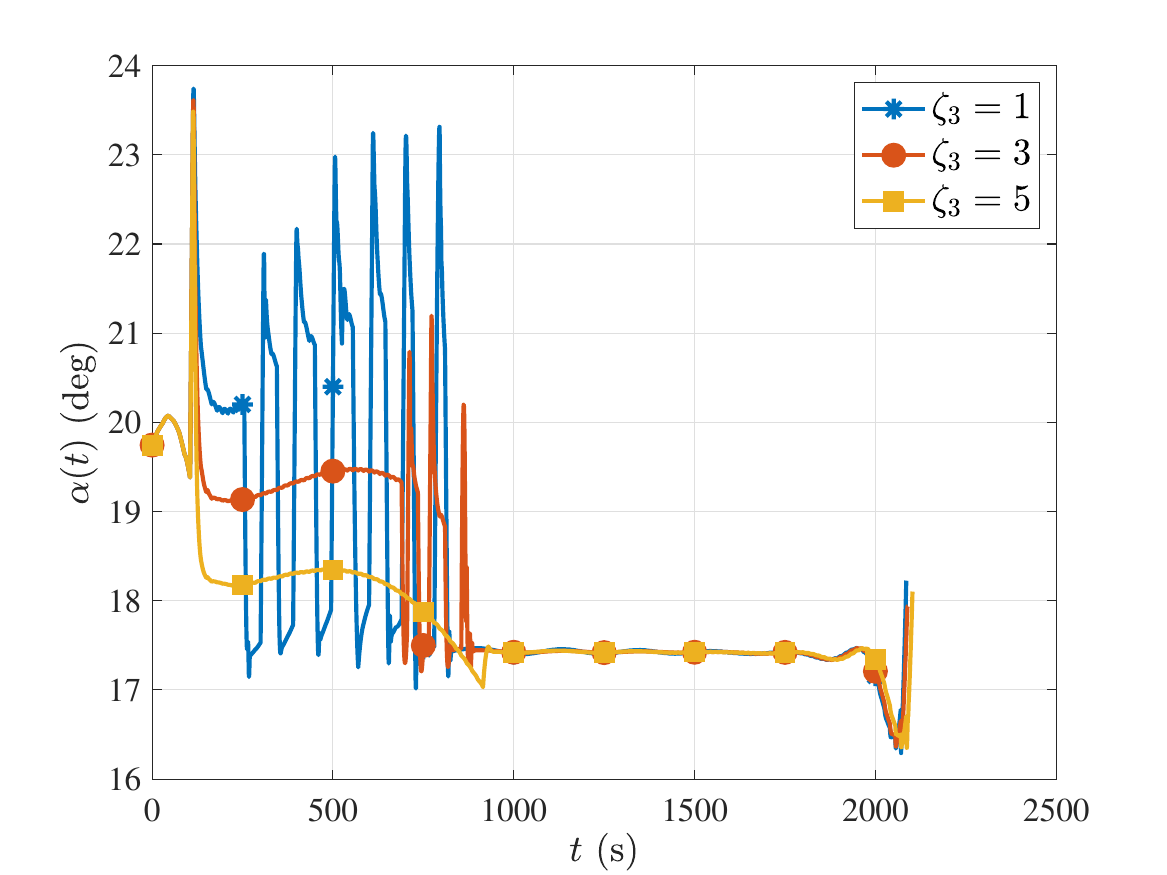}}~~~\subfloat[ Bank angle, $\sigma (t) $ vs time, $t$\label{fig:sigma phu compared}]{\includegraphics[scale=0.4]{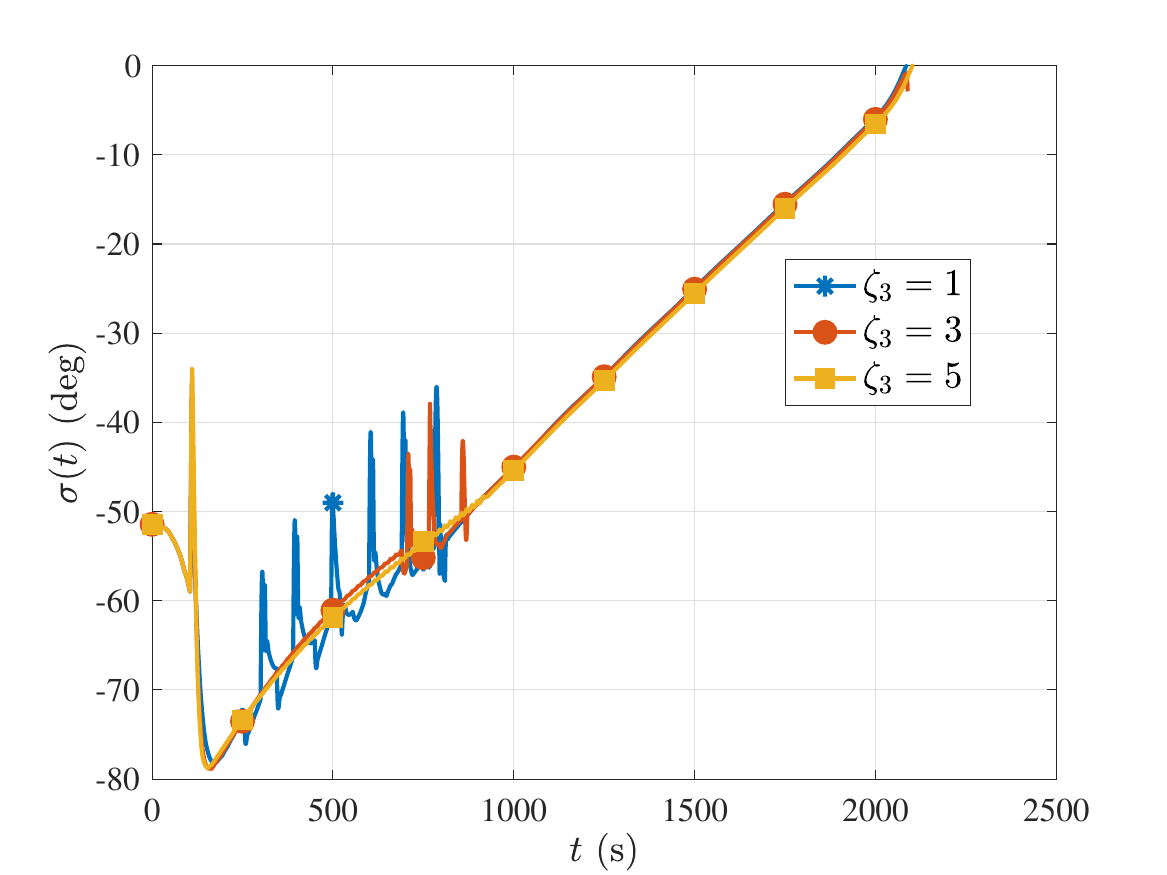}}
  \caption{Angle of attack and bank angle profiles for varied values of $\zeta_3$. }\label{fig:controls phu comp}
\end{figure}

 In Fig.~\ref{fig: states phu compared}, the trajectories obtained with the phugoid penalty have much smoother, glide-like behavior compared to the trajectories obtained using the CCG\&C method in Section~\ref{subsec:comparing weights}. In Fig.~\ref{fpa phu compared}, the phugoid penalty is shown to significantly reduce the large oscillations in flight path angle profiles.  The variation between the trajectories obtained using different values of $\zeta_3$ is seen more prominently observing the constraint profiles in Fig.~\ref{fig:path cons phu compared}. In Fig.~\ref{fig:heating rate phu compared}, as $\zeta_3$ is increased, the small oscillations in the path constraint profiles are attenuated. For the weight $\zeta_3 = 5$, the oscillations are removed entirely from the path constraint profiles.

The incorporation of the phugoid penalty reduces the drastic changes in the angle of attack and bank angle profiles that were seen in Section~\ref{subsec:comparing weights}. In Fig.~\ref{fig:controls phu comp}, as the value of $\zeta_3$ increases, the number and magnitude of the abrupt changes in the angle of attack and bank angle decrease. For the value of $\zeta_3 = 5$, the drastic changes are almost removed entirely, resulting in more realistic, executable angle of attack and bank angle profiles.

Table \ref{tab: phu comp results} provides key numerical results obtained using the CCG\&C method with different values of $\zeta_3$. As $\zeta_3$ is increased, performance increases, as the value of the final latitude $\phi(t_f)$ increases. As $\zeta_3$ is increased, the maximum heating rate also increases, but remains sufficiently below the limit such that feasibility with respect to the constraints is maintained across the entire horizon.
Finally, the errors in the terminal value of the state are all on the same order of magnitude as the results presented in Table \ref{tab: weights comp results} and within reasonable ranges.

\begin{table}[ht]
\centering
\caption{Numerical results obtained for varying weights on the phugoid penalty. \label{tab: phu comp results}}
\renewcommand{\baselinestretch}{1.25}\normalsize\normalfont
\begin{tabular}{c c c c c c}\hline\hline
  $\zeta_3$ &$\phi(t_f)$  & $\delta h (t_f)$  & $\delta v (t_f)$  & $\delta \gamma (t_f)$  & $\max (\dot{Q})$ \\
   & (deg) & (m) & (m/s) & (deg) & (MW/m$^2$) \\\hline
  $1$ & $33.382$ &$ 10.571$ & $ -0.8086$ & $ 0.0924$ & $0.8214$ \\
   $3$ & $33.700$ &$15.925$ & $ -0.9858$ & $0.1065$ & $0.8211$ \\
   $5$ & $33.898$ & $13.454$ & $ -0.9099$ & $0.1008$ & $0.8299$ \\\hline\hline
\end{tabular}
\end{table}

\subsection{Monte-Carlo Campaign}\label{subsect:Monte-Carlo simulations}

\begin{figure}[h]
\centering

\subfloat[ Altitude, $h$ vs speed, $v$ \label{fig:alt versus speed Monte}]{\includegraphics[scale=0.4]{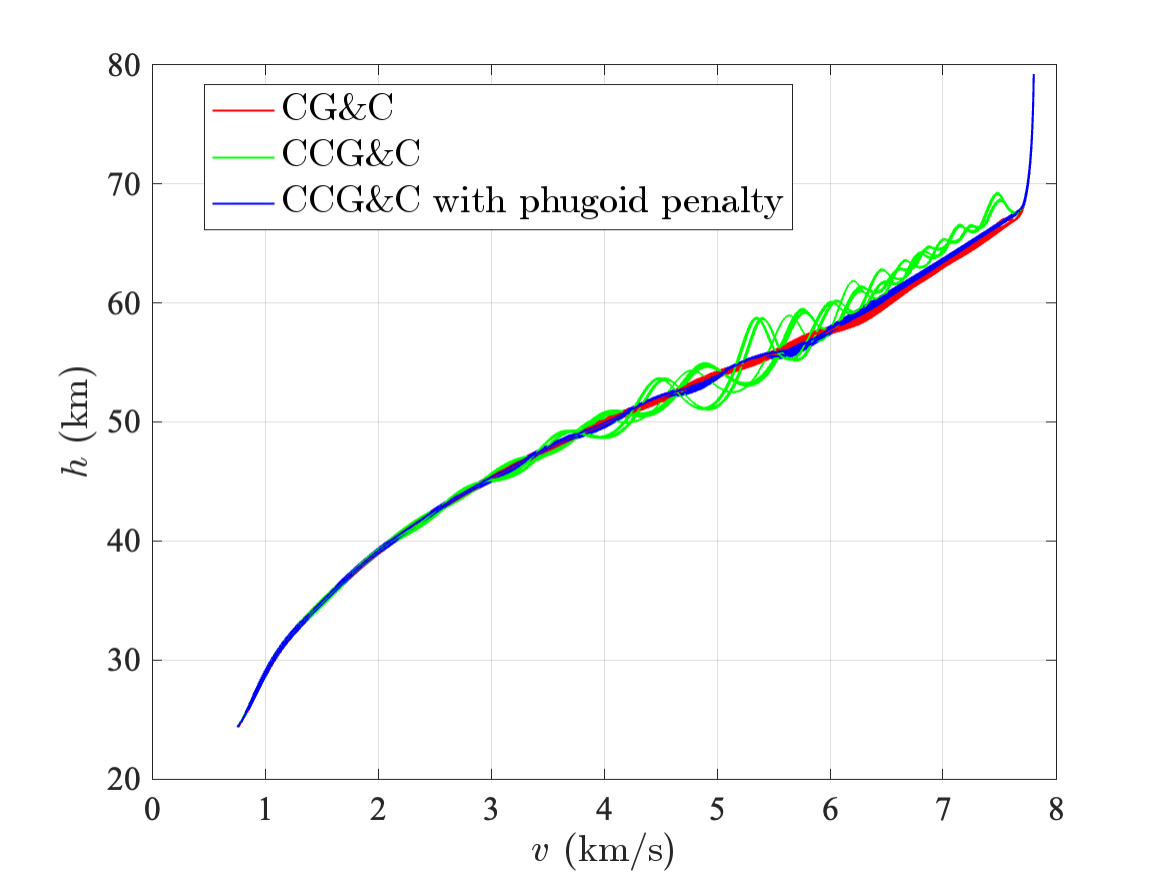}}~~~\subfloat[ Latitude, $ \phi $ vs longitude, $\theta$ \label{lat v lon Monte}]{\includegraphics[scale=0.4]{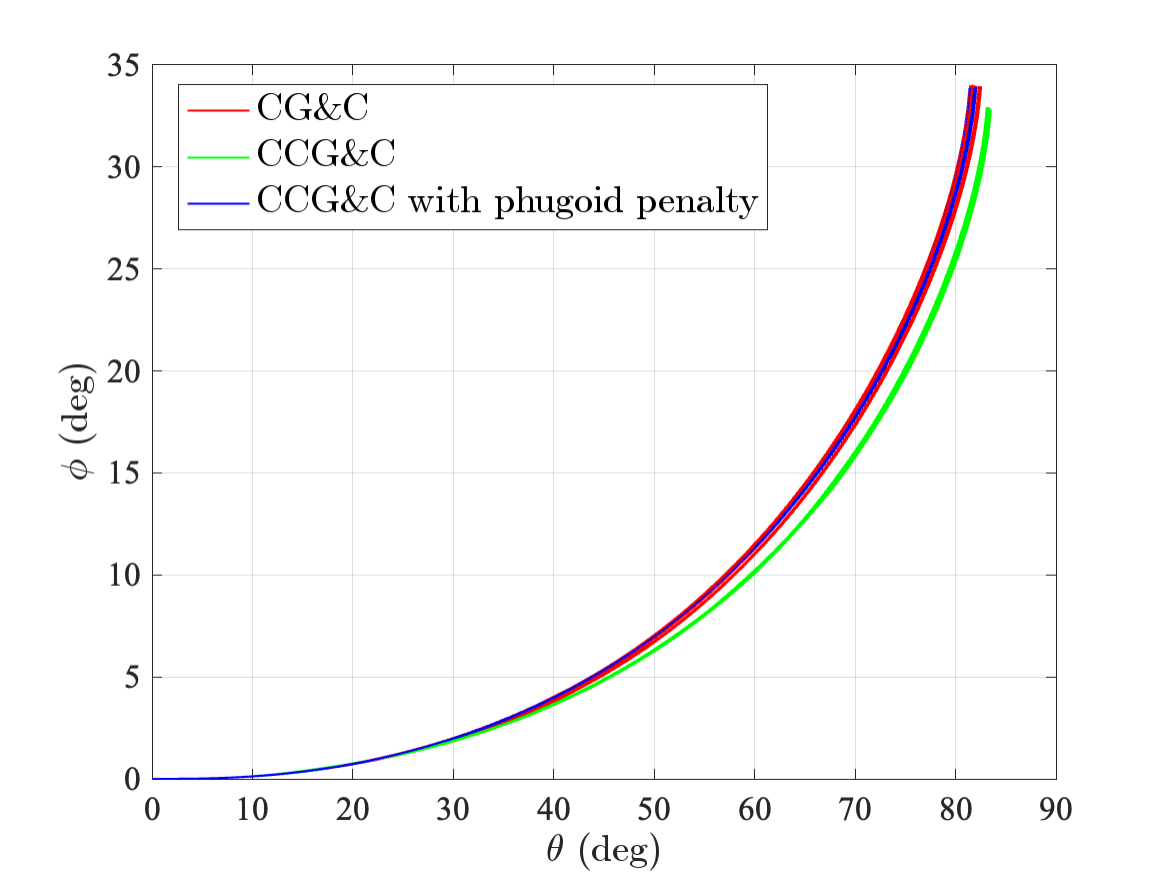}}

\subfloat[ Flight path angle, $ \gamma (t) $ vs time, $t$ \label{fpa Monte}]{\includegraphics[scale=0.4]{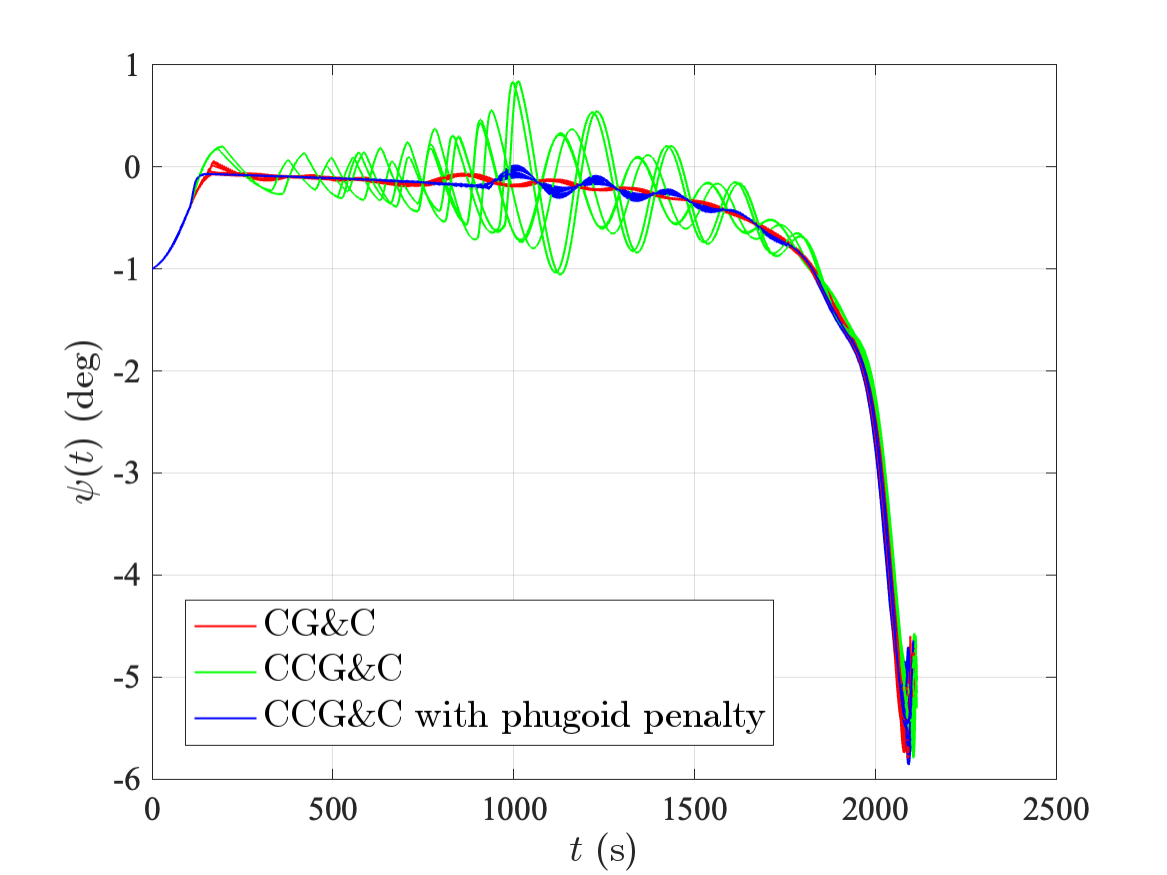}}

  \caption{Monte-Carlo campaign: state profiles obtained from simulation of perturbed dynamics. }\label{fig: states Monte}
\end{figure}

\begin{figure}[h]
\centering

\subfloat[ Heating rate, $\dot{Q}(t)$ vs time, $t$ \label{fig:heating rate weights compared full view Monte}]{\includegraphics[scale=0.4]{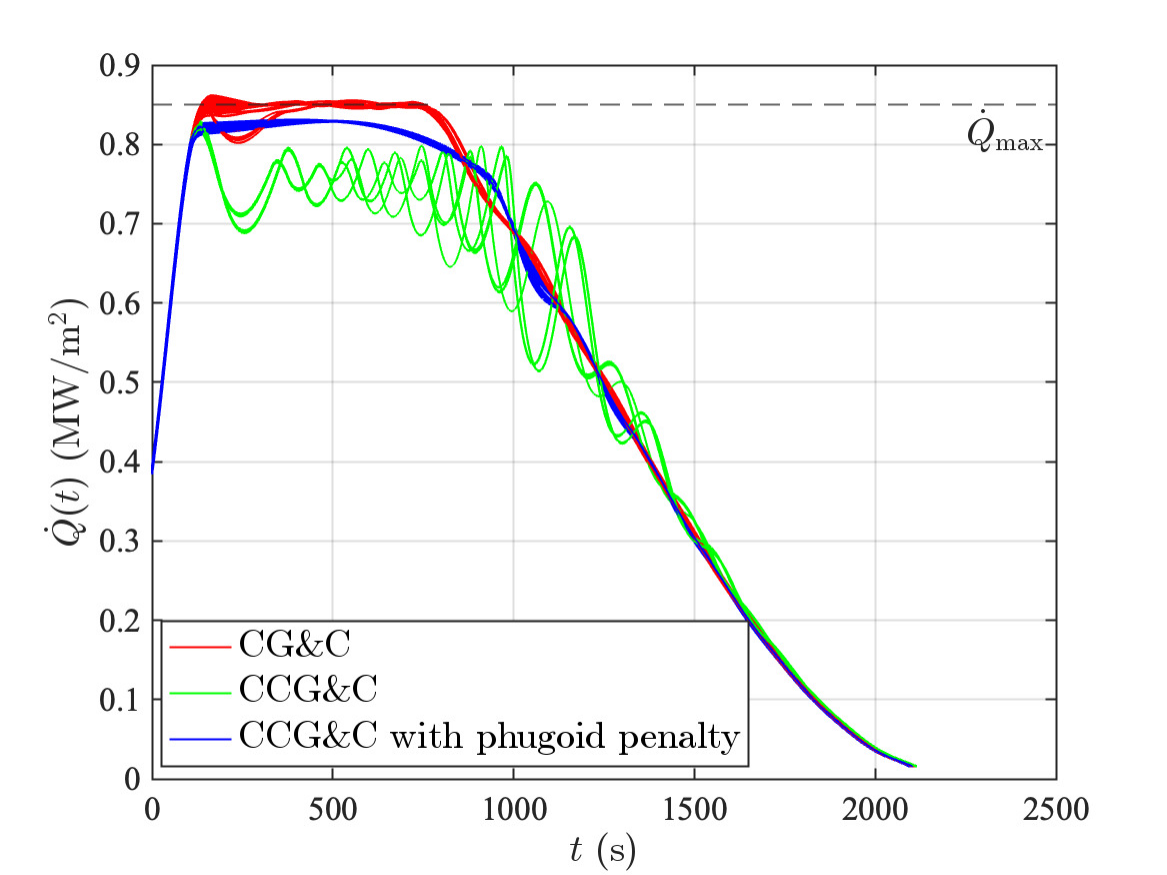}}~~~\subfloat[ Enlarged view of constraint violations\label{fig:heating rate enlarged Monte}]{\includegraphics[scale=0.4]{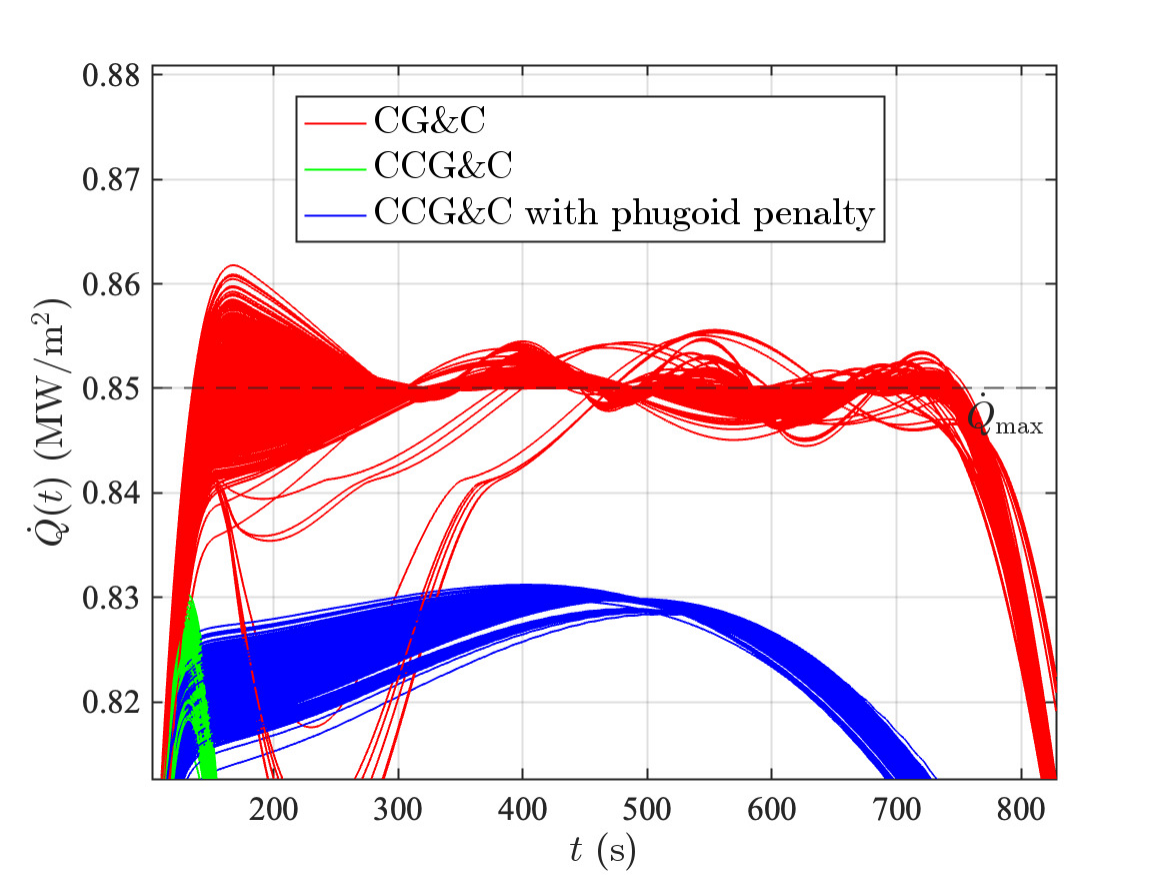}}

  \caption{Monte-Carlo campaign: heating rate profiles obtained from simulation of perturbed dynamics. }\label{fig:heating rate Monte}
\end{figure}

\begin{figure}[h]
\centering

\subfloat[ Angle of attack, $\alpha(t)$ vs time, $t$ \label{fig:aoa Monte}]{\includegraphics[scale=0.4]{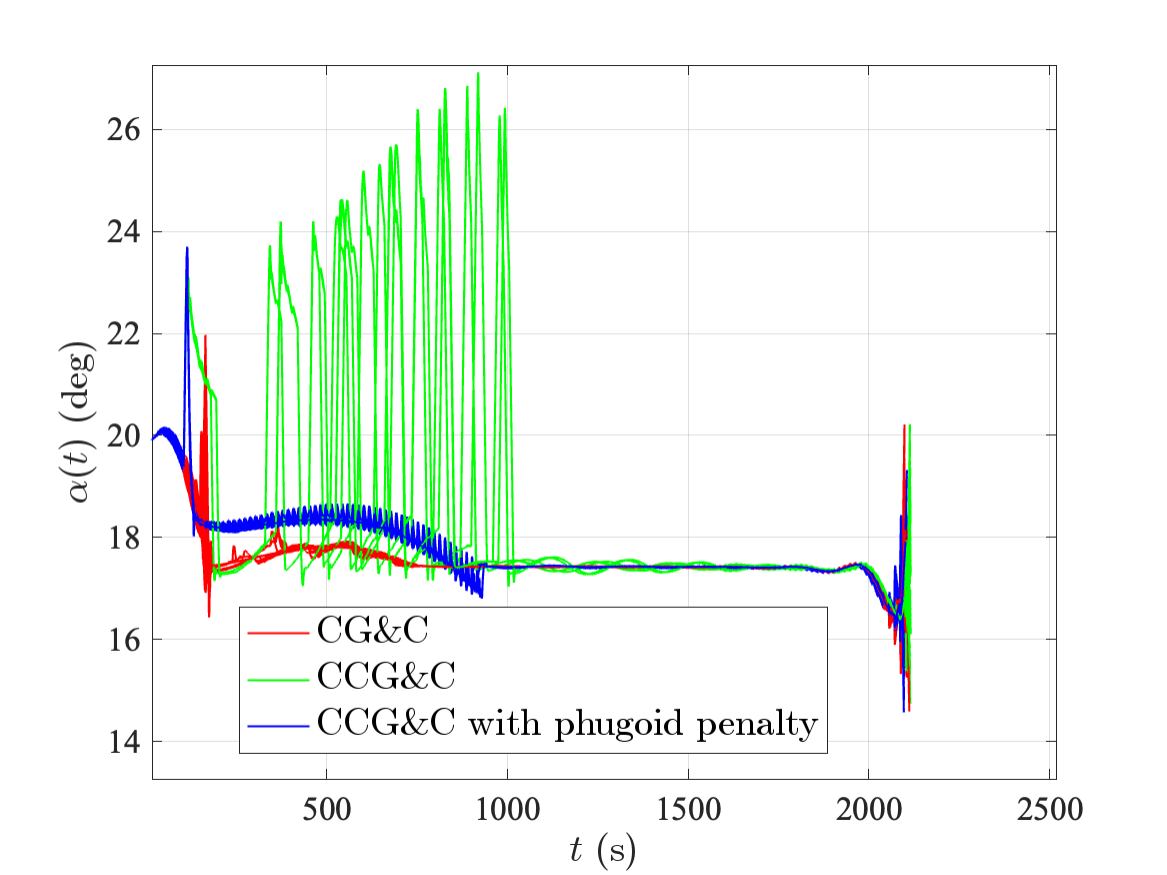}}~~~\subfloat[ Bank angle, $\sigma (t) $ vs time, $t$\label{fig:sigma Monte}]{\includegraphics[scale=0.4]{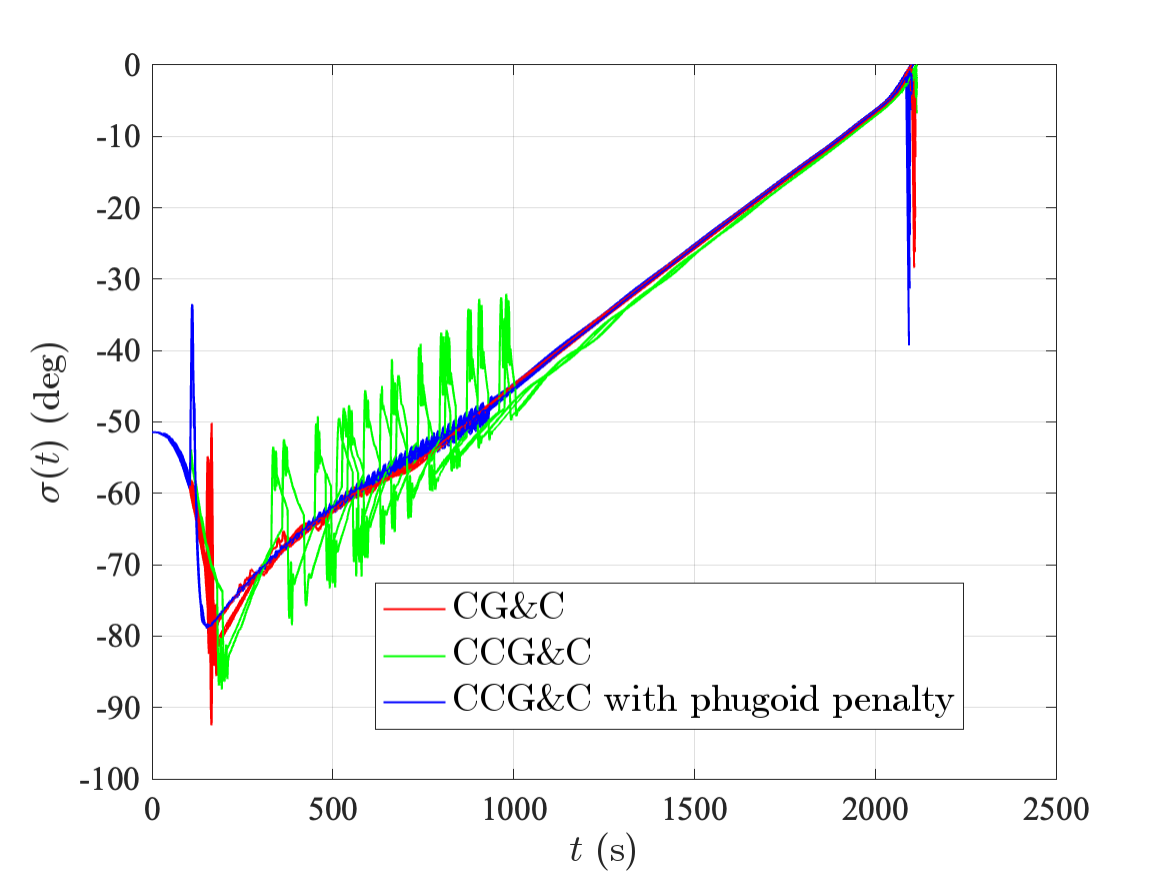}}
  \caption{Monte-Carlo campaign: angle of attack and bank angle profiles. }\label{fig:controls Monte}
\end{figure}

\begin{figure}[h]
\centering

\subfloat[CCG\&C method \label{fig:histogramfinalaltccgc}]{\includegraphics[scale=0.28]{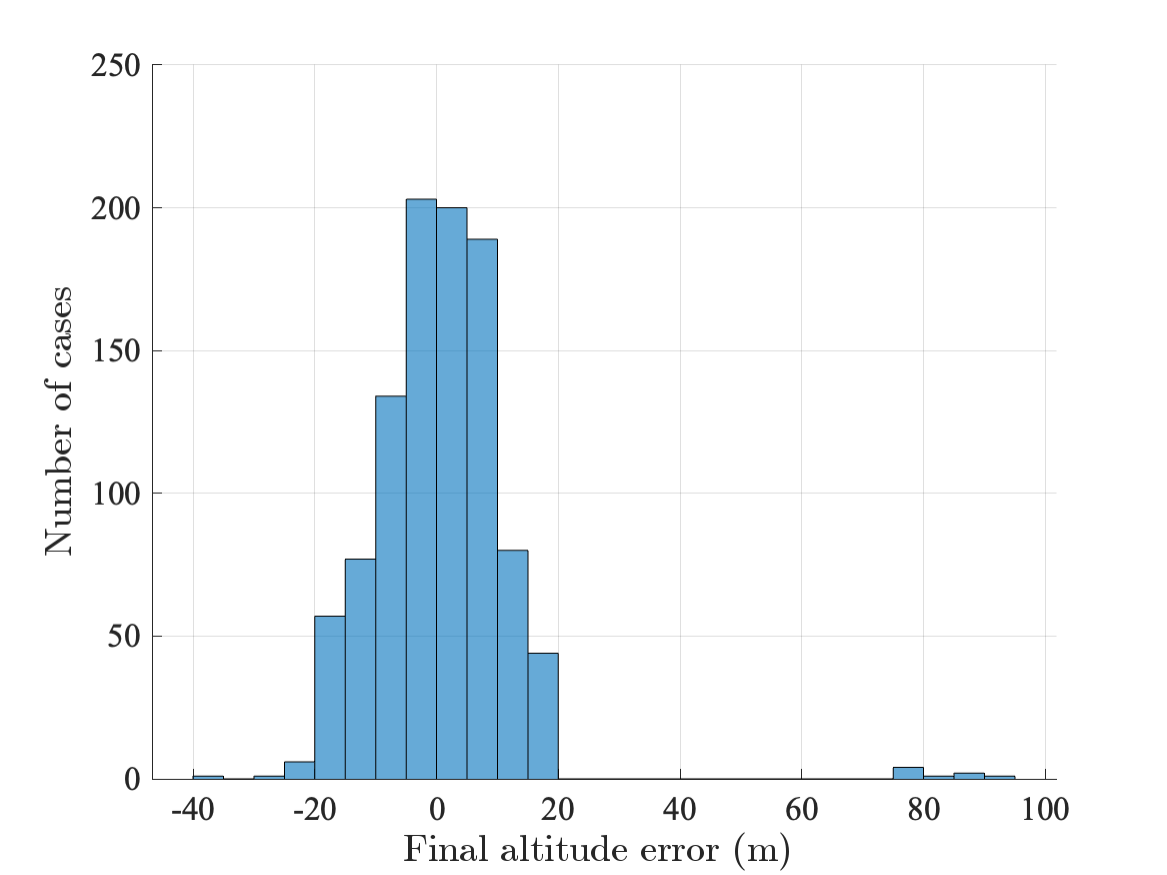}}~~~\subfloat[CCG\&C method with phugoid penalty \label{fig:histogramfinalaltccgc plus phugoid}]{\includegraphics[scale=0.28]{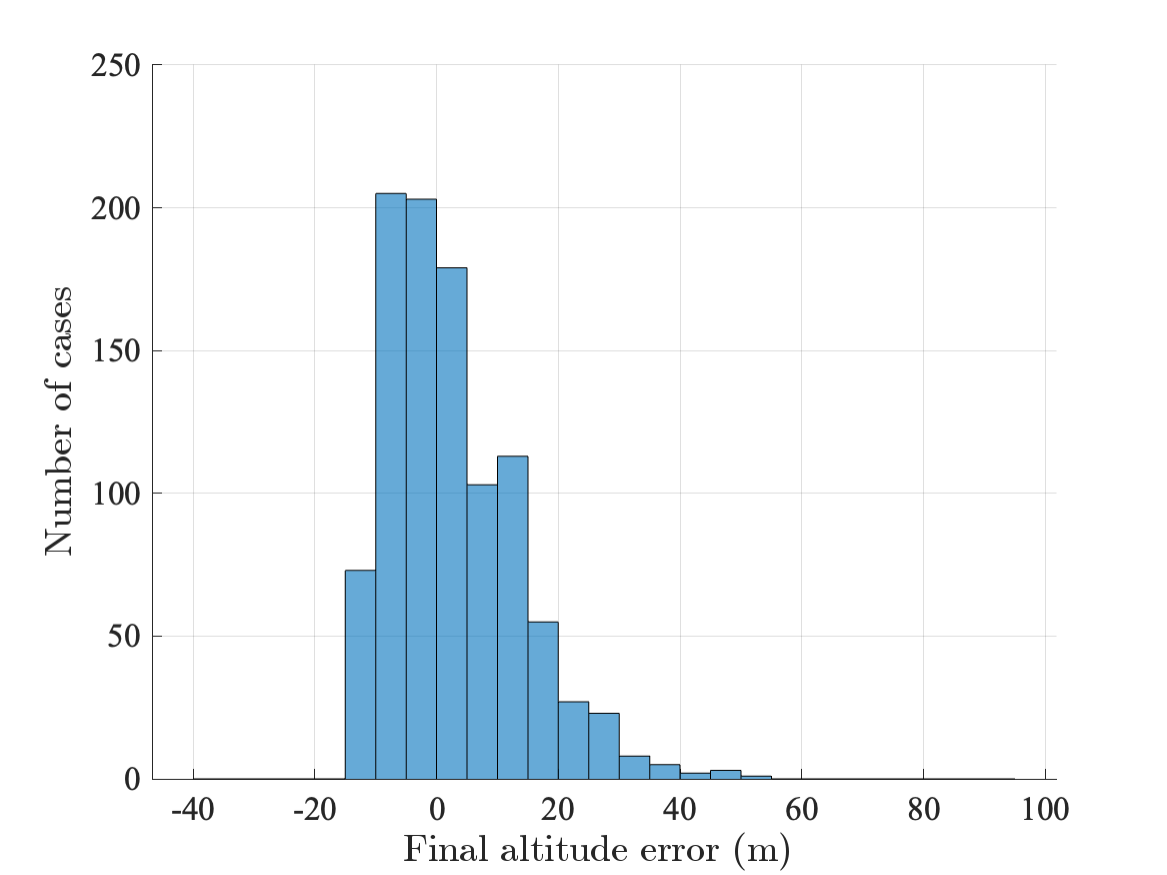}}~~~\subfloat[CG\&C method\label{fig:histogram final alt mitz}]{\includegraphics[scale=0.28]{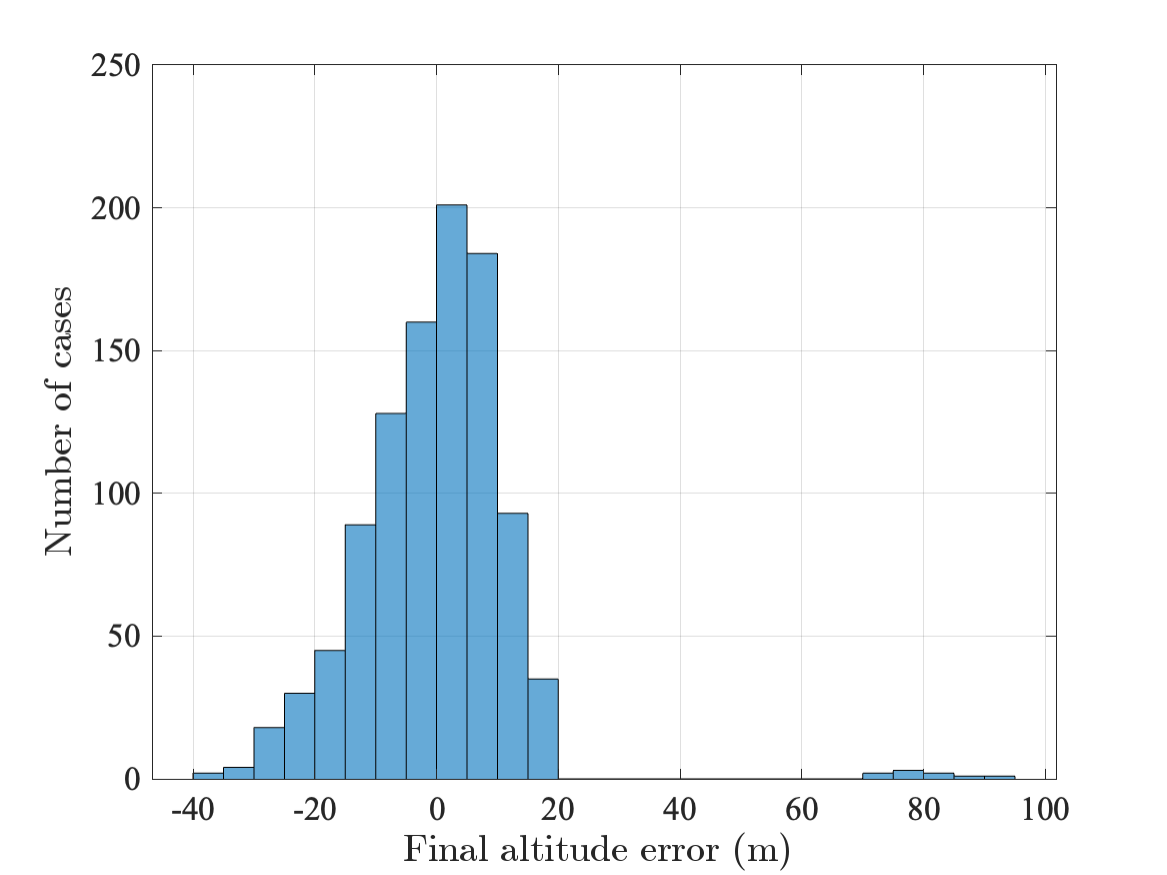}}

  \caption{Final altitude error distributions for Monte-Carlo campaign. }\label{fig:histograms final alt errors}
\end{figure}

\begin{figure}[h]
\centering

\subfloat[CCG\&C method \label{fig:histogramfinal speed ccgc}]{\includegraphics[scale=0.28]{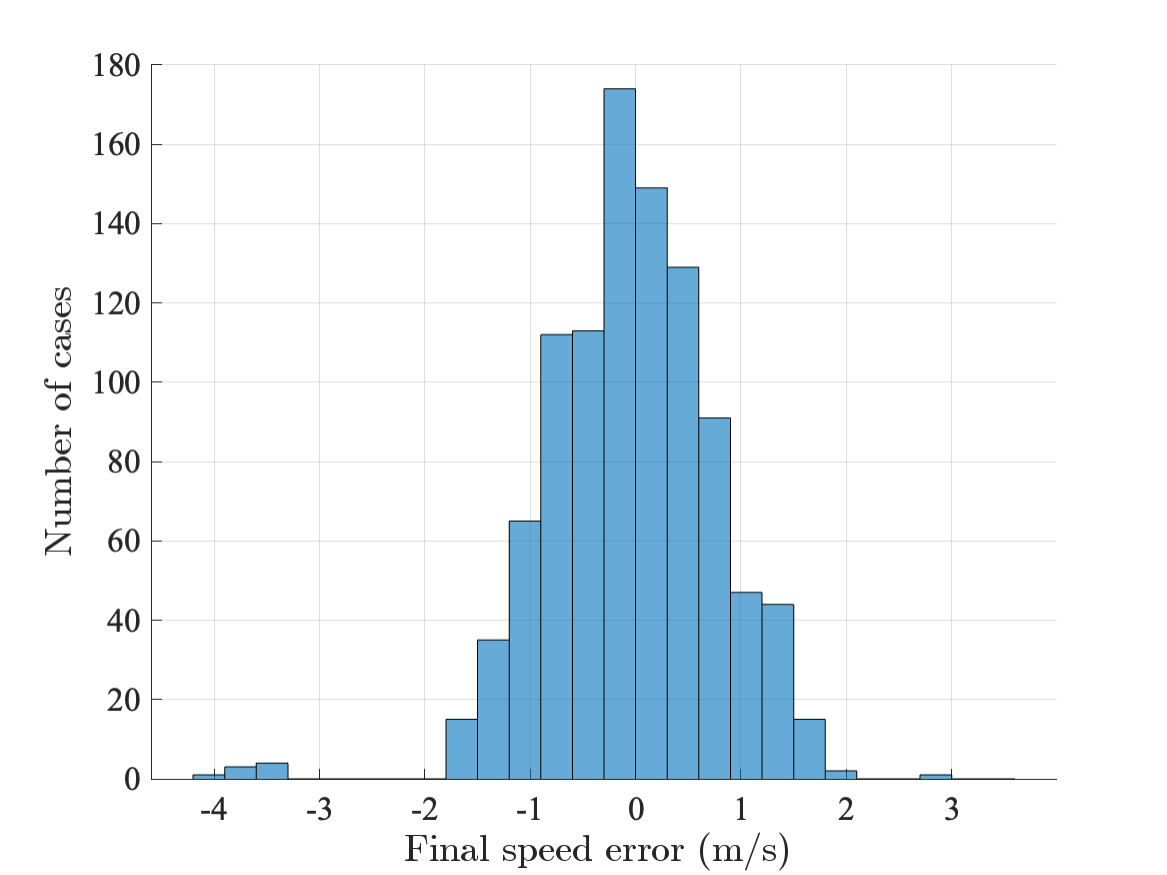}}~~~\subfloat[CCG\&C method with phugoid penalty \label{fig:histogramfinal speed ccgc with phu}]{\includegraphics[scale=0.28]{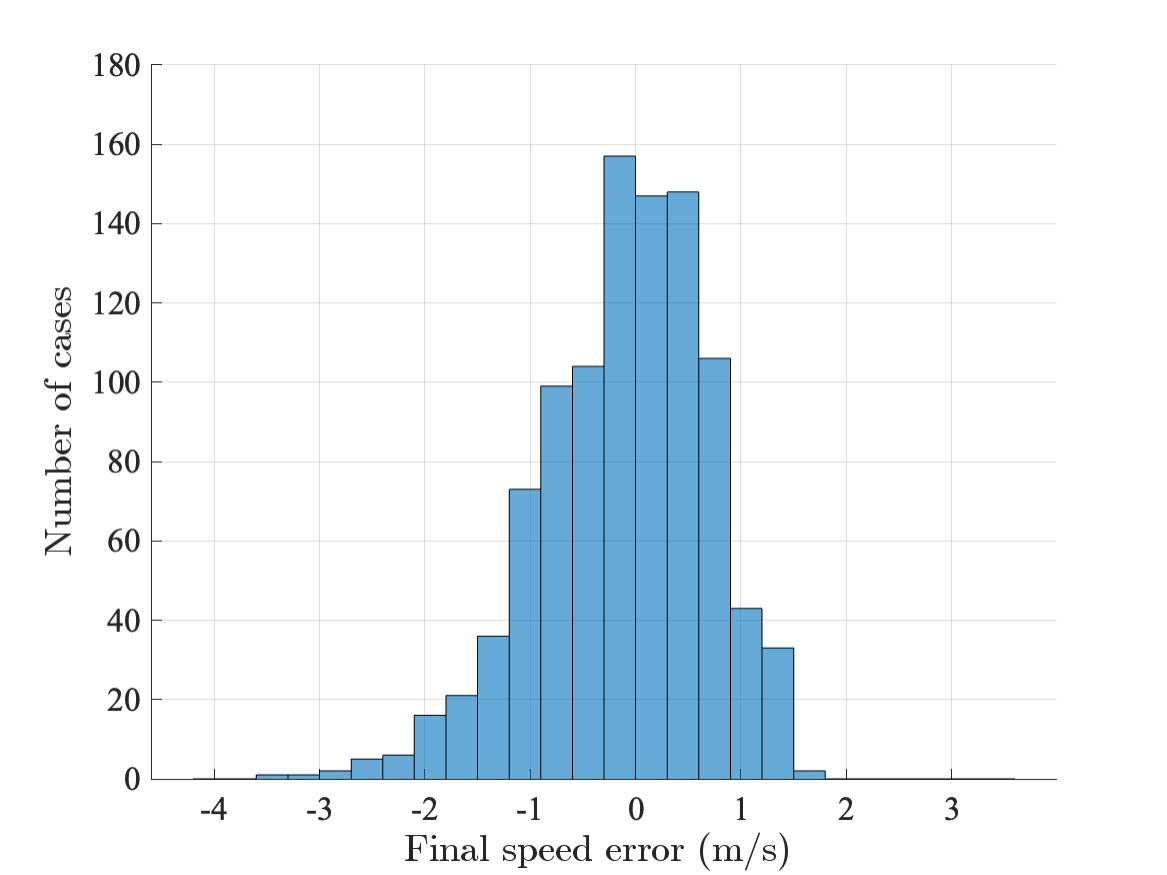}}~~~\subfloat[CG\&C method\label{fig:histogram final speed mitz}]{\includegraphics[scale=0.28]{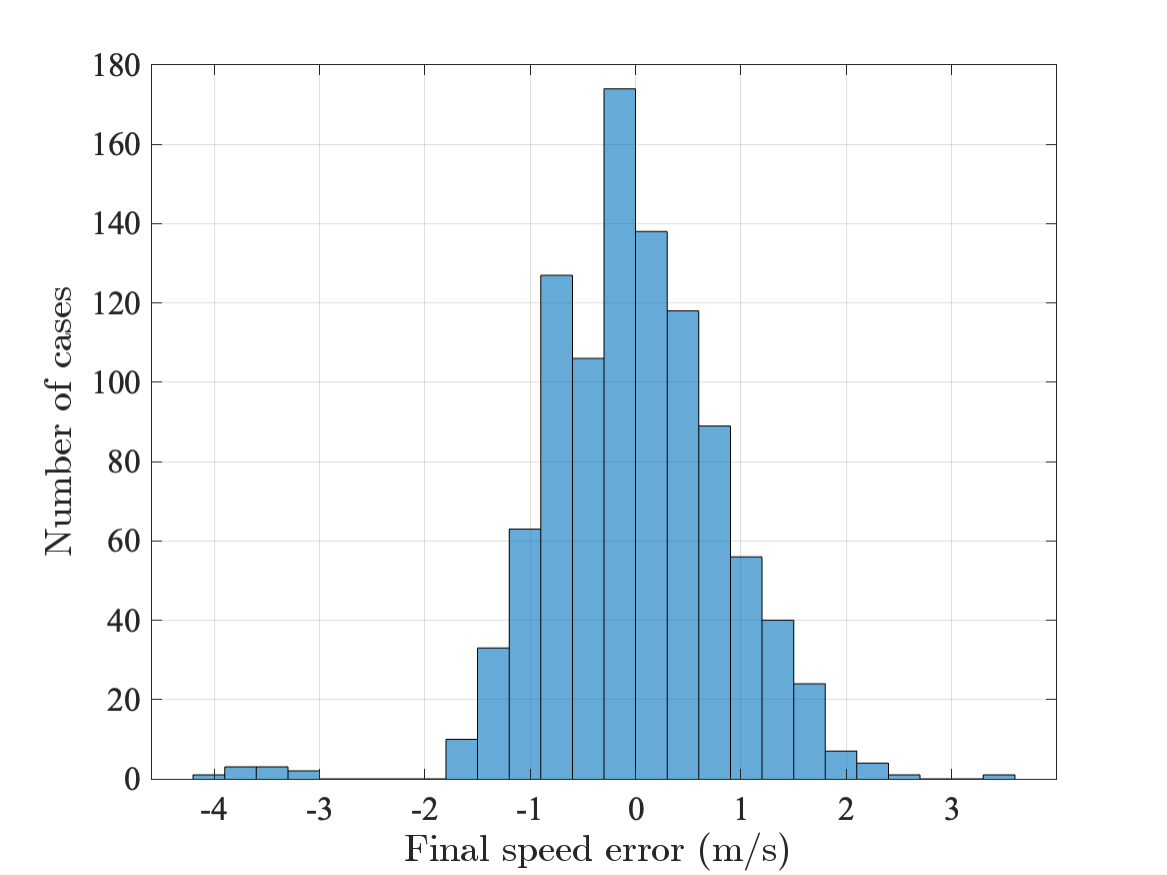}}

  \caption{Final speed error distributions for Monte-Carlo campaign. }\label{fig:histograms final speederrors}
\end{figure}

\begin{figure}[h]
\centering

\subfloat[ CCG\&C method \label{fig:histogramfinal fpa ccgc}]{\includegraphics[scale=0.28]{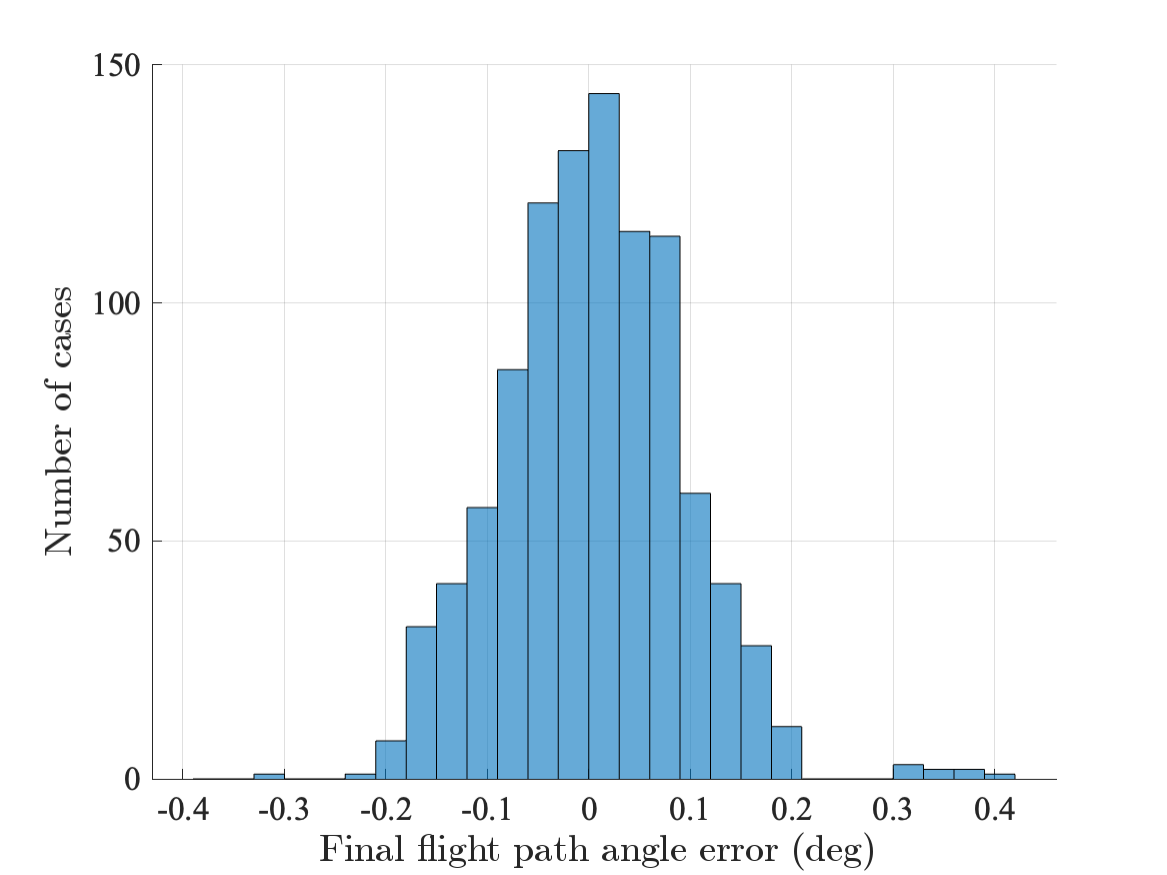}}~~~\subfloat[ CCG\&C method with phugoid penalty \label{fig:histogramfinal fpa ccgc phu}]{\includegraphics[scale=0.28]{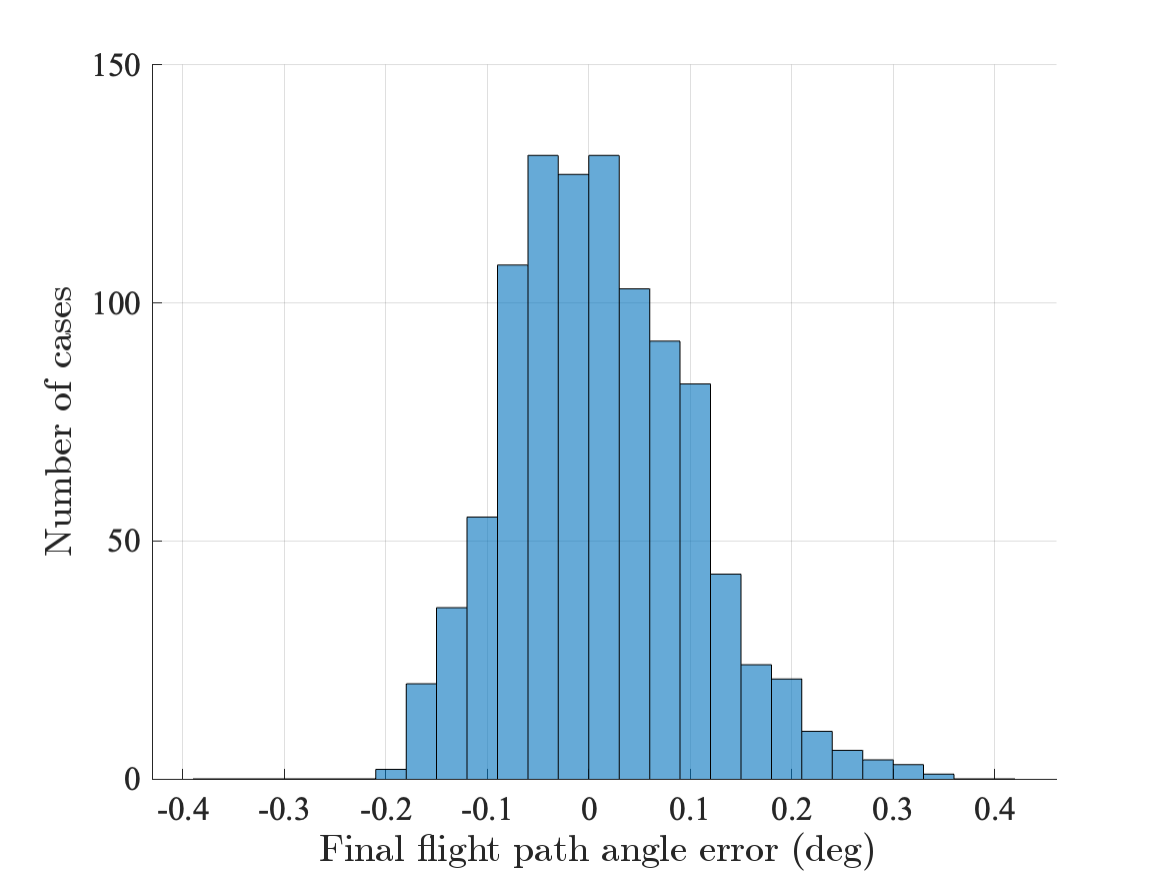}}~~~\subfloat[ CG\&C method\label{fig:histogram final fpa mitz}]{\includegraphics[scale=0.28]{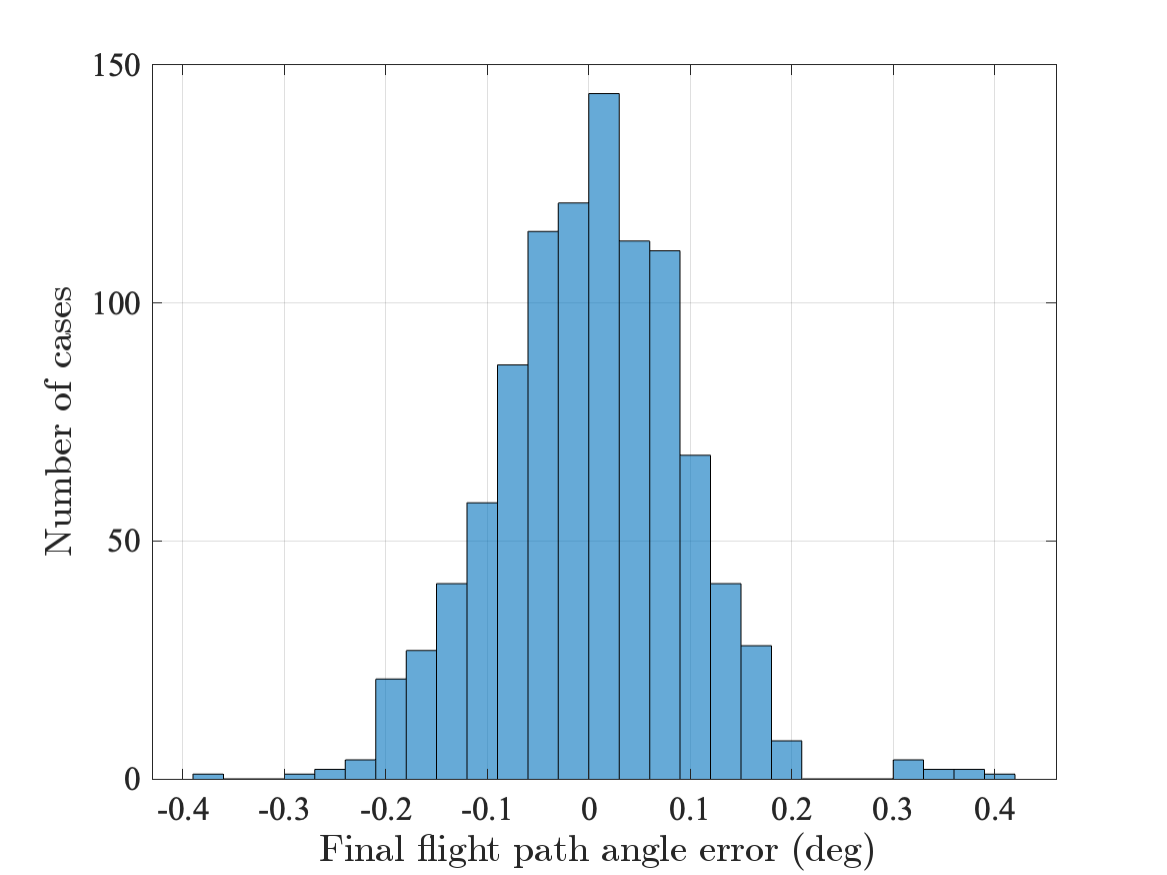}}

  \caption{Final flight path angle error distributions for Monte-Carlo campaign. }\label{fig:histograms final flight path angle errors}
\end{figure}

\begin{figure}[h]
\centering
{\includegraphics[scale=0.4]{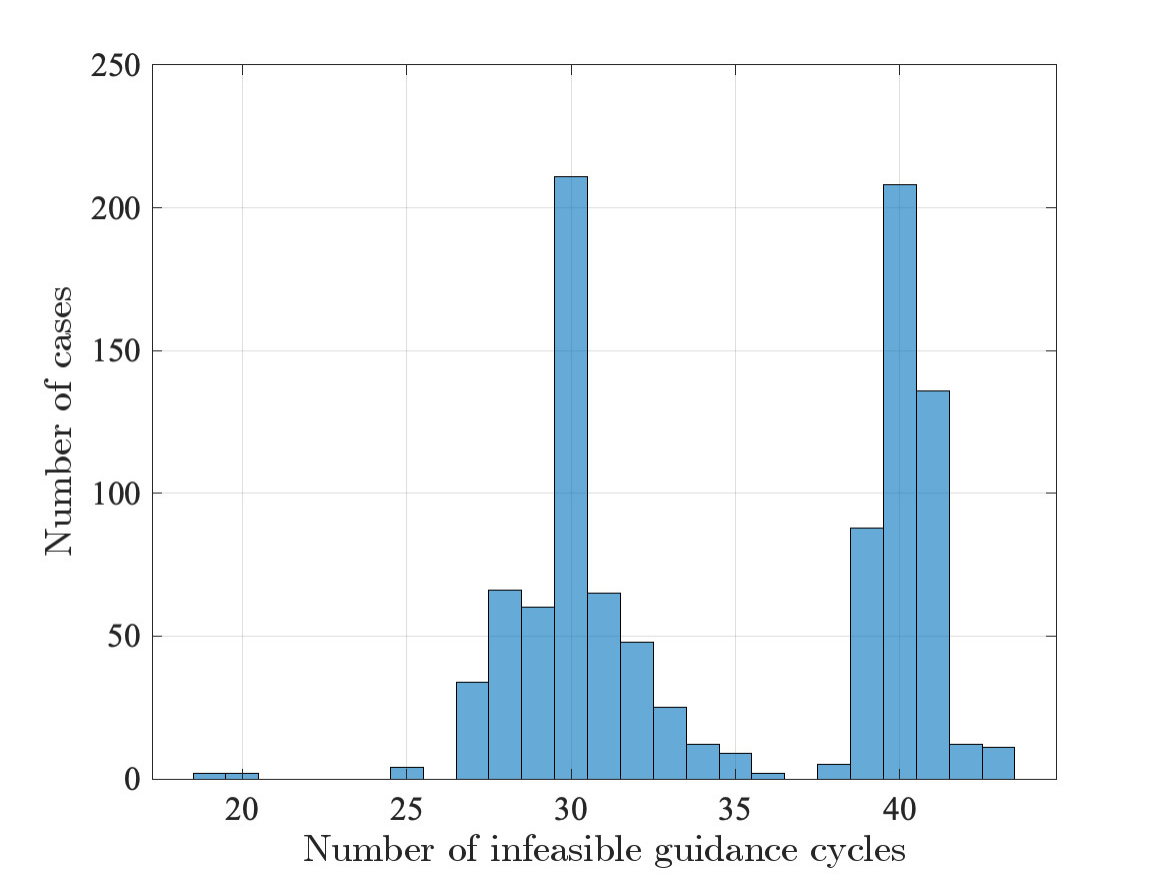}}
  \caption{Distribution of the number of infeasible guidance cycles per Monte-Carlo trial for the CG\&C method of Ref.\cite{Dennis2019}. }\label{fig: infeasible cycles}
\end{figure}

In this section, a Monte-Carlo campaign of 1000 perturbed cases of surface level atmospheric density is used to evaluate and compare the performance of the CCG\&C method developed in Section~\ref{sect:new-cgs-method} with the performance of the CG\&C method of Ref.~\cite{Dennis2019}. Results are shown for the CCG\&C method for two cases: first, using the modified objective functional as formulated in Section~\ref{sect:new-cgs-method}, and second, using the modified objective functional including the phugoid penalty given in Section~\ref{subsec:Reducing phugoids}. At the start of each Monte-Carlo trial, a perturbed value of the surface level atmospheric density was drawn from a normal distribution, where one standard deviation was chosen to be $1 \% $ of the nominal value, that is $\tilde{\rho_0} \sim \C{N}(\rho_0,(0.01\rho_0)^2)$. The perturbed value of the surface level atmospheric density was held constant for each Monte-Carlo trial. Guidance updates were made every $\Delta T = 15$ seconds. For both cases using the CCG$\&$C method, the results shown were obtained using the following parameters: $\xi = 0.9$, $\eta_i = c_{\max,i}$, $\beta = 5$, and $(\zeta_1,\zeta_2)=(0.8,0.2)$. For the case of the CCG\&C method with the phugoid penalty, the results shown were obtained using the parameters $\zeta_3=5$ and $k=3$.

Figure \ref{fig: states Monte} shows the perturbed state profiles obtained from the Monte-Carlo campaign. As expected, the trajectories obtained using the CCG\&C method without the phugoid penalty contain large oscillations induced by the abrupt changes in the angle of attack and bank angle profiles shown in Fig.~\ref{fig:controls Monte}. The incorporation of the phugoid penalty results in trajectories that are more smooth, similar to the profiles obtained using the CG\&C method. In particular, the incorporation of the phugoid penalty results in latitude versus longitude profiles nearly identical to those obtained using the CG\&C method. In fact, from Table~\ref{tab:final error results monte}, the mean final latitude for the CCG\&C method with the phugoid penalty is only slightly lower than the mean final latitude obtained using the CG\&C method, and is more than one degree higher than the mean final latitude obtained using the CCG\&C method without the phugoid penalty.

Figure \ref{fig:heating rate Monte} shows the perturbed heating rate profiles for the three guidance approaches. Observing Fig.~\ref{fig:heating rate enlarged Monte}, the CG\&C method experiences prolonged periods of constraint violation.  In fact, the trajectories obtained using the CG$\&$C method experienced heating rate constraint violations for all 1000 Monte-Carlo trials. The consequences of the heating rate constraint violations are demonstrated in Fig.~\ref{fig: infeasible cycles}, which shows the distribution of the number of infeasible guidance cycles that occurred during each Monte-Carlo trial for the CG\&C method of Ref.~\cite{Dennis2019}. The CG\&C method of Ref.\cite{Dennis2019} experienced an average of 35 infeasible guidance cycles per trial, and a maximum of 44 infeasible guidance cycles per trial. Conversely, the CCG$\&$C method (both with and without the phugoid penalty) prevented heating rate violations for all 1000 Monte-Carlo trial, ensuring feasibility with respect to the path constraints across all guidance cycles. Table~\ref{tab:final error results monte} shows that the mean maximum heating rate for the CCG\&C method without the phugoid penalty is only slightly lower than the mean maximum heating rate for the CCG\&C method with the phugoid penalty. Notably, the phugoid penalty reduced the standard deviation of the maximum heating rate value by an order of magnitude as compared to both the CCG\&C method without the phugoid penalty and the CG\&C method of Ref.~\cite{Dennis2019}.

Numerical results pertaining to the errors in the terminal value of the state are given in Table~\ref{tab:final error results monte}, and the distributions of these errors are shown in Figs.~\ref{fig:histograms final alt errors}-\ref{fig:histograms final flight path angle errors}. While the CCG\&C method without the phugoid penalty resulted in the lowest mean absolute errors, it also resulted in the largest maximum absolute errors in the terminal value of the state. The maximum absolute errors for the CCG\&C method without the phugoid oscillations are only slightly larger than the maximum absolute errors experienced by the CG\&C method. The incorporation of the phugoid penalty lowered the maximum absolute errors in the terminal value of the state, with a significant reduction in the maximum absolute error of the terminal altitude. Additionally, the CG\&C method experienced the largest mean absolute error and standard deviation in the errors, likely due to the infeasible guidance cycles where guidance updates could not be made. It is noted, however, the errors in the terminal state are all on the same order of magnitude, and within reasonable ranges.

The mean computational time required to solve the optimal control problem for each guidance cycle was $0.2133$~s, with a standard deviation of $0.6876$~s for the CCG\&C method without the phugoid penalty. For the CCG$\&$C method with the phugoid penalty, the mean computational time was $0.4361$~s, with a standard deviation of $0.8558$~s. The CG$\&$C method of Ref.~\cite{Dennis2019} had a mean computational time of $0.09038$~s with a standard deviation of $0.0377$~s. For all three guidance approaches, the mesh truncation and remapping strategy described in Ref.~\cite{Dennis2019} was employed, which results in smaller meshes and lower computational times in the later guidance cycles as the time horizon shrinks. It is important to note that on guidance cycles where the constraints were violated, the optimal control problem could not be solved. Therefore, for the CG$\&$C method, an accurate comparison of the computational time cannot be made as many of the early guidance cycles (where computational times would have been the highest) were infeasible. 

\clearpage
\begin{table}[ht]
\centering
\caption{Numerical results for Monte-Carlo Campaign\label{tab:final error results monte}}
\renewcommand{\baselinestretch}{1.25}\normalsize\normalfont
\begin{tabular}{llccc}\hline\hline
Quantity & Unit & CCG\&C method & CCG\&C method & CG\&C method  \\
&  &  & with phugoid penalty & \\\hline
  Mean $|\delta h(t_f)|$ & m & $7.861$ &$8.455$ & $8.906$   \\
  Max $|\delta h(t_f)|$ & m & $91.8435$ & $54.052$ & $91.355$ \\
  Standard deviation $\delta h(t_f)$ & m & $ 11.580$ &$   10.958$ &  $12.9985$ \\
    Mean $|\delta v(t_f)|$ & m/s & $0.6083$ & $ 0.6290$ & $  0.6467$   \\
  Max $|\delta v(t_f)|$ & m/s & $ 4.1651$ & $ 3.316$ & $4.1568$ \\
  Standard deviation $\delta v(t_f)$ &  m/s & $ 0.7920$ & $   0.7951$ &  $0.8449$\\
Mean $|\delta \gamma(t_f)|$ & deg & $ 0.0696$ & $0.0713$ & $0.0729$   \\
  Max $|\delta \gamma(t_f)|$ & deg & $0.4027$ & $0.3574$ & $ 0.4025$ \\
  Standard deviation $\delta \gamma(t_f)$ & deg & $ 0.0889$ & $ 0.0897$ &  $0.0934$ \\
 Mean $ \phi(t_f)$ & deg & $ 32.849$ & $33.895$ & $ 34.000$   \\
  Max $ \phi(t_f)$ & deg & $ 32.934$ & $  33.904$ & $34.0120$ \\
  Standard deviation $\phi(t_f)$ & deg & $0.0683$ & $ 0.00453$ &  $  0.00575$ 
  \\
     Mean $ \max(\dot{Q})$ & MW/m$^2$ & $0.8233$ & $0.8298$ & $0.8521$   \\
  Max $ \max(\dot{Q})$ & MW/m$^2$ & $ 0.8304$ & $0.8312$ & $0.8618$ \\
  Standard deviation $\max(\dot{Q})$ & MW/m$^2$ & $0.00217$ & $ 3.252\times 10^{-4}$ &  $0.00195$
  \\\hline\hline
\end{tabular}
\end{table}

\section{Discussion}\label{discussion}

The results of Section~\ref{results} highlight several key features of the constrained computational guidance and control (CCG\&C) method developed in this paper. The first key property of the method pertains to the path constraint functions. The method is designed such that, as the trajectory of the actual system approaches a constraint limit, the optimal control problem is modified to include an additional term in the objective functional. This additional term is intended to increase the margin between the constraint function and its limit. Consequently, when implementing the control obtained by solving the modified optimal control problem, the distance between the path constraint function and the path constraint limit will increase.  Figure ~\ref{fig:heating rate Monte} shows that, via appropriate choice of the weights in the terms of the modified objective functional, a significant margin is maintained between the heating rate and its limit throughout the entire trajectory.

The second key aspect of the CCG\&C method is that an accurate terminal state is attained. Specifically, the results of the Monte-Carlo campaign performed in Section~\ref{subsect:Monte-Carlo simulations} show that the mean absolute differences between the terminal conditions and their specified values (both with and without the phugoid penalty) were less than $10$ m in altitude, $1$ m/s in speed, and $0.1$ deg in flight path angle. These results demonstrate that the CCG\&C method is able to compensate for a wide range of perturbations. 

A third important aspect of the developed method relates to the control commands. It is important to note that solving the modified optimal control problem can alter the structure of the commanded inputs as compared to the commands obtained from the reference optimal control problem. Section~\ref{subsec:comparing weights} showed drastic changes in the angle of attack and bank angle profiles when solving the modified optimal control problem, which induced undesirable phugoid oscillation. Consequently, in Section~\ref{subsec:Reducing phugoids}, an additional term was incorporated in the modified objective functional, which is intended to attenuate the phugoid oscillations. Figure~\ref{fig: states Monte} shows that appropriately chosen weights on this additional term produce smooth glide-like trajectories that maintain a large margin between the constraint function and limit. The additional term improved the performance of the CCG\&C method, increasing average final latitude by more than one as compared to the CCG\&C method without the phugoid penalty.

The fourth important aspect pertains to the computational efficiency of the method. It is noted that the computational time required to solve the optimal control problem at each guidance cycle did increase slightly when additional terms were incorporated in the objective functional. While the mean computational times were lower than the guidance cycle duration, it is noted that different approaches could be taken to improve computational time required. First, executing the code in a compiled language could reduce the computational time. Second, different optimization techniques may improve speed. The results of this study were obtained using IPOPT, which employs an interior point method. Interior point methods have excellent performance in full-Newton mode, however, they perturb the initial guess away from the starting point. As a result, even with the previous solution as a good initial guess, several iterations may be required before converging to the optimal solution. It is possible other optimization techniques, such as sequential quadratic programming may have faster convergence given a good initial guess.  

Finally, it is noted that other approaches could be taken to avoid large changes in the control structure obtained from the modified optimal control problem as compared to the reference. For example, additional terms could be incorporated in the modified objective function to minimize the rates of change in the control, the magnitude of the control, or deviations from the reference control. Alternatively, more restrictive limits could be placed on the control or control rates, though this approach may overly restrict the feasible solution space.

\section{Conclusions}\label{conclusions}

A method has been described for constrained computational guidance and control in which an optimal control problem is transcribed into a nonlinear programming problem and re-solved at each guidance cycle. The method of this paper focuses on preventing constraint violation when the actual system deviates from the reference solution. Constraint violation is prevented by modifying the objective functional of the optimal control problem by adding a term that is intended to increase the margin (in the feasible direction) between the constraint function and its limit. To strike a balance between performance and constraint satisfaction, the weights in the modified objective functional are chosen adaptively at the start of each guidance cycle based on the proximity of trajectory of the actual (perturbed) system to the constraint limit. The method is demonstrated on a reusable launch vehicle entry problem where the objective is to maximize the terminal cross-range. The developed method is first studied for varying weights in the modified objective functional for a single perturbation. Then, an additional penalty term is introduced to attenuate potential phugoid oscillations. A Monte-Carlo campaign is performed to study the constrained computational guidance and control method, both with and without the phugoid penalty, and the performance is compared against a previously developed method for computational guidance and control. The results of this study demonstrate that the constrained computational guidance and control method developed in this paper is viable for path-constrained optimal guidance and control.

\section*{Acknowledgments}

The authors gratefully acknowledge support for this research from The U.S.~Air Force Research Laboratory under grant FA8651-24-1-0004, the U.S.~National Science Foundation under grant CMMI-2031213, and from the U.S.~Office of Naval Research  under grant N00014-22-1-2397.  x

\clearpage

\bibliographystyle{aiaa}

\end{document}